\documentclass[11pt]{article}
\usepackage{cancel}
\usepackage{rotating}
\usepackage{url}
\usepackage[margin=2.7 cm,nohead]{geometry}
\usepackage{hyperref} 
\usepackage[ruled]{algorithm2e} 
\usepackage{epstopdf}
\usepackage{cite}
\usepackage{algorithmic}
\usepackage{color}
\usepackage{url}
\usepackage{cite}
\usepackage{stackrel,cancel}
\usepackage{amsmath,amsthm,amsfonts,amssymb,amscd, amsxtra, mathrsfs}
\usepackage{amsfonts,amssymb}
\usepackage{amsmath}
\usepackage{graphicx,epstopdf}
\usepackage{algorithm2e}
\usepackage[caption=false]{subfig} 
\usepackage{algorithmic, enumitem}

\newtheorem{theorem}{Theorem}[section]
\newtheorem{lemma}[theorem]{Lemma}
\newtheorem{definition}[theorem]{Definition}

\newtheorem{proposition}[theorem]{Proposition}

\newtheorem{remark}[theorem]{Remark}

\usepackage{hyperref}

\begin{document}
%%%%%%%%%%%%%%%%%%%%%%%%%%%%%%%%%%%%%%%%%%%%%%%%%%%%%%%%%%%%%
\title{A subgradient method with non-monotone line search\thanks{Submitted to the editors in 2022-04-21, 11:57.}}
\author{O. P. Ferreira \thanks{Instituto de Matem\'atica e Estat\'istica, Universidade Federal de Goi\'as,  CEP 74001-970 - Goi\^ania, GO, Brazil, E-mails:  {\tt  orizon@ufg.br}},
G. N. Grapiglia \thanks{ Universit\'e Catholique de Louvain, ICTEAM/INMA, Avenue Georges Lemaître, 4-6/ L4.05.01, B-1348, Louvain-la- Neuve, Belgium, E-mail: {\tt geovani.grapiglia@uclouvain.be}},
E. M. Santos \thanks{PhD student of Instituto de Matem\'atica e Estat\'istica, Universidade Federal de Goi\'as,  CEP 74001-970 - Goi\^ania, GO, Brazil, E-mails:  {\tt  eliandersonsantos@discente.ufg.br}},
J. C. O. Souza \thanks{Department of Mathematics, Federal University of Piau\'{i}, Teresina, PI, Brazil, E-mails: {\tt joaocos.mat@ufpi.edu.br}
  }
}
\maketitle
% REQUIRED
\begin{abstract} 
In this paper we present a subgradient method with non-monotone line search for the minimization of convex functions with simple convex constraints. Different from the standard subgradient method with prefixed step sizes, the new method selects the step sizes in an adaptive way. Under mild conditions asymptotic convergence results and iteration-complexity bounds are obtained. Preliminary numerical results illustrate the relative efficiency of the proposed method.
\end{abstract}
%%%%%%%%%%%%%%%%%%%%%%%%%%%%%%%%%%
\section{Introduction}
%%%%%%%%%%%%%%%%%%%%%%%%%%%%%%%%%%
The subgradient method for solving non-differentiable convex optimization  problems  has its origin  in the 60's,  see  \cite{Ermolev1966, Shor1985}.  Over the years it has been the subject of much interest, attracting the attention  of the scientific community working on convex optimization.   One of the factors that explains the interest in the subgradient method lies in its simplicity and ease of implementation for a wide range of problems,  where the sub-differential of the objective function can be  easily computed. In addition, this method has  low  storage cost and ready exploitation of separability and sparsity, which makes it attractive in solving  large-scale problems.  For all these reasons, several variants of this method have emerged and properties of it have been discovered throughout the years, resulting in a wide literature on the subject; including  for exemple \cite{Bertsekas1999, GoffinKiwiel1999, KiwielBook1985,Kiwiel2004,nedic_bertsekas2001rate,NedicBertsekas2010,Nesterov2014, PolyakBook} and the references therein. 

The classical subgradient method employs a predefined sequence of step sizes. Standard choices include a constant step size and also sequences that converge to zero sublinearly. In this paper, we propose a subgradient method with adaptive step sizes for the minimization of convex functions with simple convex constraints in which a projection on it is easily computed. At each iteration, the selection of the step size is done by a line search in the direction opposite to the subgradient. Since, in general, this direction is not a descent direction, we endow the method with a non-monotone line search. The possible increase in the objective function
values at consecutive iterations is limited by a sequence of positive parameters that implicitly controls the step sizes. Remarkably, it is shown that the proposed method enjoys convergence and complexity properties similar to the ones of the classical subgradient method when the sequence that controls the non-monotonicity satisfies suitable conditions. Illustrative numerical results are also presented. They show that the proposed non-monotone method compares favorably with the classical subgradient method endowed with usual prefixed step sizes.

The organization of the paper is as follows. In Section \ref{sec:int.1}, we present some notation and basic results used in our presentation. In Section \ref{Alg1} we describe the subgradient method with non-monotone line search and  the main results of the present paper, including the converge theorems and iteration-complexity bounds.   Some numerical experiments are provided in Section \ref{Sec:NumExp}.   We conclude the paper with some remarks in Section \ref{Sec:Conclusions}.
 
%%%%%%%%%%%%%%%%%%%%%%%%%%%%%%%%%%%%%%%%%%%%%%%%%%%%%%%%%%
\section{Preliminaries} \label{sec:int.1}
In this section we present some notations, definitions, and results that will be used throughout the paper, which  can be found  in \cite{Amir, Lemarechal}.

{\it Denotes ${\mathbb N}:=\{1, 2, 3,  \ldots \}$.} A function $f\colon\mathbb{R}^{n}\to \mathbb{R} $ is said to be $\sigma$-strongly convex with modulus $\sigma\geq 0$ if  $f(\tau x + (1-\tau)y)\leq \tau f(x) + (1-\tau) f(y)-\frac{\sigma}{2}\tau(1-\tau)\|x-y\|^2$, for all $x,y\in {\mathbb{R}^{n}}$ and $\tau \in [0,1]$.  For $\sigma = 0$ we  say that $f$ is a convex function.

\begin{proposition}\label{pr:secv}  
The function  $f:\mathbb{R}^{n}\to \mathbb{R} $ is  $\sigma$-strongly convex  with modulus $\sigma\geq 0$ if and only if $f(y)\geq f(x) + \langle v, y-x \rangle + ({\sigma}/{2}) \| y-x\Vert ^{2}$, for all $x,y\in \mathbb{R} ^{n} $ and all $v\in \partial f(x)$.
\end{proposition}
 A function $f:\mathbb{R}^{n} \to  \mathbb{R}$ is {\it $L_{f,{\cal C}}$-Lipschitz continuous on ${\cal C}\subset \mathbb{R}^{n}$ if  there exist a constant $L_{f,{\cal C}}>0$}  such that  $|f(x)-f(y)|\leq L_{f,{\cal C}}\|x-y\|$,  for all $x, y \in {\cal C}$. Whenever $C=\mathbb{R}^{n}$ we set $L_{f}\equiv L_{f,\mathbb{R}^{n}}$.
\begin{proposition} \label{pr:cdd}
Let $f:\mathbb{R}^{n}\to\mathbb{R}$ be a  convex. Then, for all $x\in \mathbb{R}^{n}$ the set  $\partial f(x)$ is a non-empty, convex, compact subset of  $\mathbb{R}^{n}$. In addition,  $f$  is $L_{f,{\cal C}}$-Lipschitz function on ${\cal C}\subset \mathbb{R}^{n}$ if and only if  $\|v\|\leq L_{f,{\cal C}}$ for all $v\in \partial f(x)$  and $x\in  {\cal C}$.
\end{proposition}
\begin{remark}
In view of Proposition~\ref{pr:cdd}, if ${\cal C}\subset \mathbb{R}^{n}$ is a compact set then $f$ is a $L_{f,\cal C}$-Lipschitz function on ${\cal C}\subset \mathbb{R}^{n}$
for some $L_{f,\cal C}>0$.
\end{remark} 
\begin{definition} 
Let ${\cal C} \subset {\mathbb R}^n$ be a closed convex set. The {\it projection map}, denoted by ${\cal P}_{\cal C}: {\mathbb R}^n \rightrightarrows {\cal C}$,  is  defined as follows ${\cal P}_{\cal C}(y):=\arg\min \{\|y-z\|:~z\in {\cal C}\}.$
\end{definition}
The next lemma  presents an important property  of the projection.
%%%%%%%%%%%%%%%%%%
\begin{proposition}  \label{pr:nexp}
	Let  $y \in {\mathbb R}^n$ and  $z \in {\cal C}$.  Then, we have $ \left\|{\cal P}_{\cal C}(y) - z\right\|^2 \leq \|y-z\|^2.$
\end{proposition}
\begin{definition}\label{def:fejer}
Let $S$ be a nonempty subset of $\mathbb{R}^n$. A sequence $(v^k)_{k\in\mathbb{N}}\subset \mathbb{R}^n$ is said to be quasi-Fej\'er convergent to $S$, if and only if, for all $v\in S$ there exists ${\bar k}\ge 0$ and a summable sequence $(\epsilon_k)_{k\in\mathbb{N}}$, such that $\|v^{k+1}-v\|^2 \le \|v^k - v\|^2+\epsilon_k$ for all $k\ge {\bar k}$. 
\end{definition}
In the following lemma, we state the main properties of quasi-Fej\'er sequences that we will need; a comprehensive study on this topic can be found in \cite{Combettes2001}.
\begin{lemma}\label{le:fejer}
Let  $(v^k)_{k\in\mathbb{N}}$  be   quasi-Fej\'er convergent to $S$.  Then, the following conditions hold:
\begin{itemize}
\item[(i)] the sequence $(v^k)_{k\in\mathbb{N}}$ is bounded;
\item[(ii)] if a cluster point ${\bar v}$ of $(v^k)_{k\in\mathbb{N}}$  belongs to $S$, then  $(v^k)_{k\in\mathbb{N}}$ converges to ${\bar v}$.
\end{itemize}
\end{lemma}
%%%%%%%%%%%%%%%%%%%%%%%%%%%%%%%%%%%%%%%%%%%%%%%%%%%%%%%%%%%%%
\section{Subgradient method with non-monotone line search} \label{Alg1}
We are interested in the following constrained optimization problem
\begin{equation}\label{eq:pr_linear}
\begin{array}{c}
\min f(x)\\
\mbox{s.t. } x\in {\cal C}, 
\end{array}
\end{equation}
where $f:\mathbb{R}^n \to \mathbb{R}$ is a convex function  and $  {\cal C}\subset \mathbb{R}^{n}$ is a closed and convex set.  {\it Denote by $\Omega^*$ the optimal set  of the problem~\eqref{eq:pr_linear} and by  $f^*$ the optimal value}.  Throughout  the paper we will consider problem~\eqref{eq:pr_linear} under the following  two assumptions:
\begin{enumerate}
\item [(H1)] $f:\mathbb{R}^{n}\to \mathbb{R}$ is  a convex function and $L_{f, {\cal C}}$-Lipschitz continuous;
\item [(H2)] $f^{*}:=\inf _{x\in {\cal C}} f(x) >-\infty .$
\end{enumerate}
 We propose the following conceptual algorithm to find  a solution  of problem~\eqref{eq:pr_linear}.\\

%\hrule
\begin{algorithm} {\bf SubGrad projection method with non-monotone line search} \label{Alg1s}
\begin{footnotesize}
\begin{description}
	\item[ Step 0.]  Fix $c>0$, $(\gamma^k)_{k\in\mathbb{N}}\subset \mathbb{R}_{++}$ a non-increasing  sequence,  $\rho, \beta\in (0,1)$ and $\alpha>0$. Choose an initial point $x_1\in C$.  Set  $\alpha_1=\alpha$ and $k=1$;
	\item[ Step 1.]  Choose $s_{k}\in\partial f(x_{k})$.  If $s_k=0$,  then STOP and returns $x_{k}$; 
	\item[ Step 2.]  Compute 
	\begin{equation} \label{eq:jk} 
\hspace{-10 pt}\ell_{k}:=\min \left\{\ell \in {\mathbb N}:~ \beta^{\ell}\alpha_{k}\leq c\beta \gamma_k , ~f( {\cal P}_{\cal C}(x_{k}-\beta^{\ell} \alpha _{k}s_{k}))\leq f(x_{k})-{\rho} \big(\beta^{\ell}{\alpha_{k} }\big)\| s_{k}\| ^{2} +\gamma_{k} \right\};  
\end{equation}
	\item[ Step 3.] Set   $x_{k+1}:={\cal P}_{\cal C}(x_{k}-\beta^{\ell_{k}}\alpha _{k}s_{k})$, $\alpha _{k+1}:= \beta^{\ell_{k}-1}\alpha_{k}$. Update   $k \leftarrow k+1$ and go to Step~1.
\end{description}
%\hrule
\end{footnotesize}
\end{algorithm}
\noindent
\vspace{0.3cm}
%%%%%%%%%%%%%%%%%%%%%%%%%%%%%%%%
\begin{remark}
It follows from \cite[Theorem 4.2.3]{Lemarechal1993} that the set where convex functions fail to be differentiable is of zero measure.  Consequently, almost every opposite direction of a subgradient is a descent direction. Therefore, we expect Algorithm~\ref{Alg1s} to be able to skip non-differentiability points that are not minimum points and then behave similarly to the gradient method with non-monotonic line search at differentiability points.    It is worth  to noting that the idea of using general non-monotone line searches in differentiable optimization, generalizing the  non-monotone searches  proposed in \cite{Grippo1986, ZhangHager2004},  have appeared  in \cite{GrapigliaSachs2017, SachsSachs2011}. A modified  version of the subgradient method with  the  non-monotone line search proposed in \cite{Grippo1986}  was considered in \cite{ krejic2022spectral2022, LoretoCrema2015, LoretoKotva2019}.\end{remark}

In the  following  lemmas   we establish   general  inequalities that are important  in our analysis. We begin presenting   the well definition of $\ell_k$ defined in {Step 2} of  Algorithm~\ref{Alg1s} and two inequalities that follows as a consequence.
\begin{lemma}\label{le:wd}  
There exists $\ell_k$ satisfying \eqref{eq:jk}. As a consequence,  the following inequalities  hold:
\begin{equation} \label{eq:ilsjk} 
\alpha _{k+1}\leq c \gamma_k, \qquad f(x_{k+1})\leq f(x_{k})-{\rho}\beta \alpha_{k+1}\| s_{k}\| ^{2} +\gamma_{k},  \quad \forall k \in \mathbb{N}, 
\end{equation}
and $x_{k+1}\in {\cal C}$,  for all $k \in \mathbb{N}$.
\end{lemma}
\begin{proof} 
Since $f$ and  the projection  ${\cal P}_{\cal C}$  are    continuous functions  and the point $x_k\in {\cal C}$, we have $\lim_{\alpha \to 0^+}(f({\cal P}_{\cal C}(x_{k}-\alpha s_{k})) - f(x_{k})+\rho \alpha \|s_{k}\|^{2})=0$. Hence, due to 
$\gamma _{k}>0$,  there exists $\eta _{k}>0$ such that $f({\cal P}_{\cal C}(x_{k}-\alpha s_{k})) - f(x_{k})+\rho \alpha \|s_{k}\|^{2}< \gamma _{k}$,  for all $\alpha \in (0,{\eta _{k}}]$,  or equivalently,
 \begin{equation}\label{eq:dscboosted}
f({\cal P}_{\cal C}(x_{k}-\alpha s_{k}))  \leq f(x_{k})-\rho \alpha \|s_{k}\|^{2}  +\gamma _{k}, \qquad \alpha \in (0, \eta_k].
\end{equation}
Hence, due to $\beta\in (0,1)$ we have $\lim_{\ell \to 0^+} \beta^{\ell}{\alpha_{k}}=0$, and since $\eta_k>0$,  we obtain that there exists  ${\ell_*}\in {\mathbb N}$ such that $\ell\geq {\ell_*}$ implies $\beta^{\ell}{\alpha_{k}}\in (0, \eta_k]$.  Therefore,  due to  \eqref{eq:dscboosted} be  hold for all $\alpha \in (0, \eta_k]$,  there exists  $\ell_{k}$ satisfying  \eqref{eq:jk}, which proves the first statement. The inequalities  in \eqref{eq:ilsjk} and inclusion $x_{k+1}\in {\cal C}$   follow from  the definitions of $x_{k+1}$ and $\alpha_{k+1}$  in Step 3. 
\end{proof}
 {\it From now on    $(x_{k})_{k\in \mathbb{N}}$ denotes the sequence  generated by Algorithm~\ref{Alg1s}}.  In the  next lemma we   recall    a classical inequality used  in  the study of subgradient methods. We give the proof here for  the sake of completeness. 
\begin{lemma}\label{Le:xkArm} 
For any   $x\in \mathbb{R}^{n}$ there  holds  
\begin{equation} \label{eq:sfi}
2\beta\alpha _{k+1}(f(x_{k})-f(x))\leq \|x_{k}-x\|^{2}-\|x_{k+1}-x\|^{2} +\beta^2 \alpha _{k+1}^{2}\|s_{k}\|^{2}, \quad \forall k \in \mathbb{N}.
\end{equation}
In addition, if  $f$ is a $\sigma$-strongly convex function  then    there  holds  
\begin{equation} \label{eq:pscsfi}
2\beta\alpha _{k+1}(f(x_{k})-f(x))\leq (1-\sigma \beta\alpha _{k+1})  \|x_{k}-x\|^{2}-\|x_{k+1}-x\|^{2}+\beta^2 \alpha _{k+1}^{2}\|s_{k}\|^{2}, \quad \forall k \in \mathbb{N}.
\end{equation}
\end{lemma}
\begin{proof}
Since the inequality \eqref{eq:pscsfi} becomes \eqref{eq:sfi} for $\sigma=0$,  it is sufficient to prove \eqref{eq:pscsfi}. 
It follows from the  definition of  $x_{k+1}$ in Step 3   of Algorithm~\ref{Alg1s},  Proposition~\ref{pr:nexp} and also definitions of $\alpha_{k+1}$ that 
\begin{align} 
\|x_{k+1}-x\|^{2}  &=\|{\cal P}_{\cal C}(x_{k}-\beta^{\ell_{k}}\alpha _{k}s_{k}) - x\|^2 \notag\\
                           &\leq \| x_{k}-\beta^{\ell_{k}}\alpha _{k}s_{k} -x\|^2. \label{eq:pclaw}\\
                           &= \|x_{k}-x\|^{2}+2\beta \alpha _{k+1}\langle s_{k},x-x_{k} \rangle +\beta^2\alpha _{k+1}^{2}\|s_{k}\|^{2}.\notag
\end{align}
Therefore, considering that  $f$ is a   $\sigma$-strongly convex function, it follows from Proposition~\ref{pr:secv}  that  $\langle s_{k},x-x_{k} \rangle\leq f(x)-f(x_k)-(\sigma/2)\|x_{k}-x\|^2$, which substituting into  \eqref{eq:pclaw} yields  \eqref{eq:pscsfi}.  The inequality \eqref{eq:sfi} follows from \eqref{eq:pscsfi} by letting $\sigma=0$.
\end{proof}
Next we present an important relationship between  $(\alpha _{k})_{k\in \mathbb{N}}$ and $(\gamma _{k})_{k\in \mathbb{N}}$. 
\begin{lemma}\label{le:bapha} 
The following  inequality holds:
\begin{equation} \label{eq:bapha}
\alpha _{k}\geq \min \left\{ \alpha _{1},\frac{\gamma_{k}}{(1+\rho)L_{f, {\cal C}}^{2}} \right\}, \quad \forall k \in \mathbb{N}.
\end{equation} 
\end{lemma}
\begin{proof}
The  inequality \eqref{eq:bapha} immediately holds for $k=1$.  Suppose by an absurd that there exists  $k\in \mathbb{N}$ such that 
\begin{equation} \label{eq:iarag}
\alpha _{k+1}< \min \left\{ \alpha _{1},\frac{\gamma_{k+1}}{(1+\rho)L_{f, {\cal C}}^{2}} \right\}.
\end{equation} 
Since we are supposing  that $(\gamma^k)_{k\in\mathbb{N}}$ is a non-increasing  sequence, using  the definition of $\alpha_{k+1}$  in Step 3   of Algorithm~\ref{Alg1s}  together with \eqref{eq:iarag},  we conclude that 
\begin{equation} \label{eq:iarag2}
\beta^{\ell_{k}-1}\alpha _{k}=\alpha _{k+1}< \min \Big\{ \alpha _{1},\frac{\gamma_{k+1}}{(1+\rho)L_{f, {\cal C}}^{2}} \Big\}\leq \min \Big\{ \alpha _{1},\frac{\gamma_{k}}{(1+\rho)L_{f, {\cal C}}^{2}} \Big\}\leq \frac{\gamma_{k}}{(1+\rho)L_{f, {\cal C}}^{2}}.
\end{equation} 
Considering that $f$ is   $L_{f,{\cal C}}$-Lipschitz continuous and  $x_k\in {\cal C}$, using Proposition~\ref{pr:nexp}  we have 
\begin{align*} 
f\big({\cal P}_{\cal C}(x_k- \beta^{\ell_{k}-1}\alpha _{k} s_k)\big)-f(x_k)&\leq L_{f,{\cal C}} \|{\cal P}_{\cal C}(x_k- \beta^{\ell_{k}-1}\alpha _{k} s_k)-x_k\|\\
                                                                                         &\leq L_{f,{\cal C}} \|x_k- \beta^{\ell_{k}-1}\alpha _{k} s_k-x_k\|\\
                                                                                         &= L_{f,{\cal C}} \beta^{\ell_{k}-1}\alpha _{k} \|s_k\|.
\end{align*} 
Using again that  $L_{f,{\cal C}}$-Lipschitz continuous, it follows from  Proposition~\ref{pr:cdd}  that   $\|s_k\|\leq L_{f,{\cal C}}$. Thus, after some algebraic manipulations, the two previous inequalities  imply that 
\begin{align*}
f\big({\cal P}_{\cal C}(x_k- \beta^{\ell_{k}-1}\alpha _{k} s_k)\big)-f(x_k) +\rho  \beta^{\ell_{k}-1}\alpha _{k}  \|s_k\|^{2} &\leq L_{f,{\cal C}}  \beta^{\ell_{k}-1}\alpha _{k}  \|s_k\|  +  \rho  \beta^{\ell_{k}-1}\alpha _{k}  \|s_k\|^{2}\\
                                                                                                                                               &\leq \ \beta^{\ell_{k}-1}\alpha _{k}  (1+\rho)L_{f, {\cal C}}^{2}. 
\end{align*}
Hence, using \eqref{eq:iarag2} we obtain that   $f({\cal P}_{\cal C}(x_k- \beta^{\ell_{k}-1}\alpha _{k} s_k))-f(x_k) +\rho  \beta^{\ell_{k}-1}\alpha _{k}  \|s_k\|^{2} <\gamma_k$, or equivalently
$$
f\big({\cal P}_{\cal C}(x_k- \beta^{\ell_{k}-1}\alpha _{k} s_k)\big)< f(x_k) -\rho  \beta^{\ell_{k}-1}\alpha _{k}  \|s_k\|^{2} + \gamma_k, 
$$
 which contradicts the definition of $\ell_k$ in \eqref{eq:jk}. Thus, \eqref{eq:bapha} holds for all $k$ and the proof is complete. 
\end{proof}
In the following we combine the inequalities \eqref{eq:ilsjk} in Lemma~\ref{le:wd} with those in Lemmas~\ref{Le:xkArm} and \ref{le:bapha} to provide an  inequality that will allow us to prove the convergence of $(x_{k})_{k\in \mathbb{N}}$ and    obtain some iteration-complexity bounds. For that,  it is convenient to define the following  positive constants for $\rho>1/2$:
\begin{equation} \label{eq:cfcb}
 {\Theta}:=\min \left\{ 1,\frac{1}{(1+\rho)L_{f, {\cal C}}^{2}} \right\}, \qquad \quad {\Gamma}:= {\Theta}\left(2\beta-\frac{\beta}{\rho}\right).
\end{equation} 
\begin{lemma}\label{le:auxineq} 
Assume that     $\Omega^*\neq \varnothing$.  Let $(x_{k})_{k\in \mathbb{N}}$ be generated by Algorithm~\ref{Alg1s} and  $x^*\in \Omega^*$.Then,  following  inequality holds:
\begin{equation}  \label{eq:sfiaux}
\Gamma \gamma_{k+1}(f(x_{k})-f^*)\leq \|x_{k}-x^*\|^{2} -\|x_{k+1}-x^*\|^{2}+  \frac{1}{\rho}\beta c\gamma_k^2,  \quad \forall k \in \mathbb{N}.
\end{equation}
In addition, if  $f$ is a $\sigma$-strongly convex function  then    there  holds  
\begin{equation}  \label{eq:psfiaux}
\Gamma \gamma_{k+1}(f(x_{k})-f^*)\leq (1-\sigma \beta \Theta \gamma_{k+1})\|x_{k}-x^*\|^{2} -\|x_{k+1}-x^*\|^{2}+  \frac{1}{\rho}\beta c\gamma_k^2,  \quad \forall k \in \mathbb{N}.
\end{equation}
\end{lemma}
\begin{proof}
First of all, note that the inequality \eqref{eq:psfiaux} becomes \eqref{eq:sfiaux} for $\sigma=0$. Then,   it is sufficient to prove the inequality \eqref{eq:psfiaux}.  It follows from Lemma~\ref{le:wd} that  $\beta \alpha_{k+1} \| s_{k}\| ^{2} \leq (f(x_{k})-f( x_{k+1})+\gamma_{k})/{\rho}$, which  combined with  inequality  \eqref{eq:pscsfi}  in Lemma~\ref{Le:xkArm}   yields  
\begin{equation} \label{eq:pcscsip}
\Big(2\beta-\frac{\beta}{\rho}\Big)\alpha _{k+1}(f(x_{k})-f^*)\leq (1-\sigma \beta\alpha _{k+1})  \|x_{k}-x^*\|^{2}-\|x_{k+1}-x^*\|^{2}+  \frac{1}{\rho}\beta\alpha_{k+1}\gamma_k.
\end{equation}
On the other hand, by using Lemma~\ref{le:bapha}, considering that   $(\gamma^k)_{k\in\mathbb{N}}$ is a non-increasing  sequence  and also using the first equality in  \eqref{eq:cfcb}  we obtain that 
\begin{equation} \label{eq:pccsiiag}
\alpha _{k+1}\geq \min \Big\{ 1,\frac{1}{(1+\rho)L_{f,{\cal C}}^{2}} \Big\} \gamma_{k+1}={\Theta}  \gamma_{k+1}.
\end{equation}
 Besides, we  know from  Lemma~\ref{le:wd}   that   $\alpha_{k+1}\leq c\gamma_k$, which  combined with  \eqref{eq:pcscsip}  and  \eqref{eq:pccsiiag} yield
\begin{equation*}
\Big(2\beta-\frac{\beta}{\rho}\Big){\Theta}  \gamma_{k+1}(f(x_{k})-f^*)\leq ({1-\sigma \beta {\Theta}\gamma _{k+1}})\|x_{k}-x^*\|^{2}- \|x_{k+1}-x^*\|^{2}+  \frac{1}{\rho}{\beta c}{\gamma_k^2},
\end{equation*}
Therefore,   taking into account \eqref{eq:cfcb},  the last inequality implies \eqref{eq:psfiaux} and the proof is concluded.
\end{proof}
\begin{remark}\label{remarkbfcnm} 
It is worth to compare the classical  inequalities \eqref{eq:sfi} and \eqref{eq:pscsfi} in Lemma~\ref{Le:xkArm}  with,  respectively,   the inequalities \eqref{eq:sfiaux} and  \eqref{eq:psfiaux} in Lemma~\ref{le:auxineq}.  This   comparison  shows  that the latter inequalities allow transfer to the sequence of non-monotonicity parameters $(\gamma _{k})_{k\in \mathbb{N}}$ the classical conditions usually  imposed on the sequence of step sizes $(\alpha _{k})_{k\in \mathbb{N}}$ that control  the behavior of  $(x _{k})_{k\in \mathbb{N}}$, see for example \cite{Amir, Bertsekas1999}. This way, the method itself will select the step sizes $\alpha_k$, which are usually prefixed in the classical formulations of sugbgradient method. In fact, for each prefixed non-increasing exogenous sequence  $(\gamma _{k})_{k\in \mathbb{N}}$, it follows from  Lemma~\ref{le:wd}, Lemma~\ref{le:bapha} and   first equality in \eqref{eq:cfcb} that Algorithm~\ref{Alg1s},  by performing a non-monotone line search,  select   the step sizes $\alpha_k$ satisfying  the following inequalities 
\begin{equation}\label{eq:bfcnm} 
 {\Theta}\gamma_{k+1}\leq \alpha_{k+1}\leq  c\gamma_{k}, \qquad k\in {\mathbb N}, 
\end{equation}
which shows that  our method is different from the  ones   that appeared in \cite{ krejic2022spectral2022, LoretoCrema2015, LoretoKotva2019}.  Moreover,  our line search  allows different choices for the sequence  $(\gamma^k)_{k\in\mathbb{N}}$  that controls the non-monotonicity.
\end{remark}
%%%%%%%%%%%%%%%%%%%%%%%%%%%%%%%%
\subsection{Convergence analysis} %\label{eq:convsec}
In this section we analyze  the behavior of the sequence  $(x _{k})_{k\in \mathbb{N}}$ under assumptions  (H1), (H2) and more two additional  assumptions. {\it The additional assumptions  will be used separately  and only when explicitly stated}. The new assumptions are as follows:
\begin{enumerate}
\item[(H3)] The sequence of non-monotonicity parameters $(\gamma^k)_{k\in\mathbb{N}}$ satisfies 
\begin{equation*} 
\lim_{N\to +\infty} \displaystyle \frac{\sum_{k=1}^{N}  \gamma_{k}^2}{\sum_{k=1}^{N}  \gamma_{k+1}}=0.
\end{equation*} 
\item[(H4)] The sequence of non-monotonicity parameters $(\gamma^k)_{k\in\mathbb{N}}$ satisfies 
\begin{equation*} 
\lim_{N\to +\infty} \displaystyle \frac{\sum_{k=1}^{N}  \gamma_{k}^2}{N \gamma_{N+1}}=0.
\end{equation*} 
\end{enumerate}
\begin{theorem} \label{th:cbn}
Assume that $\Omega^*\neq \varnothing$. Let $(x_{k})_{k\in \mathbb{N}}$ be generated by Algorithm~\ref{Alg1s} with $\rho>1/2$ and $x^*\in \Omega^*$. Then, for each fixed  $N\in {\mathbb N}$, the following inequality  hold: 
\begin{equation} \label{eq:icb2}
\min \big\{f(x_k)-f^*: ~k=1, \ldots, N\big\} \leq \frac{1}{\Gamma} \Big(\|x_1-x^*\|^2+\beta \rho^{-1}c\sum_{k=1}^N\gamma_k^2\Big)\frac{1}{\sum_{k=1}^N\gamma_{k+1}}.
\end{equation} 
Consequently,  if (H3)  holds then $\lim_{N\to +\infty}\min \big\{f(x_k)-f^*: ~k=1, \ldots, N\big\}=0$.
\end{theorem} 
\begin{proof}
Let $k\leq N$. Using the  inequality \eqref{eq:sfiaux} in  Lemma~\ref{le:auxineq} and   taking into account that  $\min \big\{f(x_k)-f^*: ~k=1, \ldots, N\big\}  \sum_{k=1}^N \gamma_{k+1}\leq   \sum_{k=1}^N   \gamma_{k+1}(f(x_{k})-f^*)$,  we obtain  that 
$$
\Gamma\min \big\{f(x_k)-f^*: ~k=1, \ldots, N\big\}  \sum_{k=1}^N \gamma_{k+1}\leq  \|x_{1}-x^*\|^{2} +    \frac{1}{\rho}\beta c \sum_{k=1}^N\gamma_k^2,
$$
which implies  \eqref{eq:icb2}. For  concluding  the proof, first note that assumption (H3)  implies that   $\lim_{N\to +\infty} {\sum_{k=1}^{N}  \gamma_{k+1}}=+\infty$. Thus, using   \eqref{eq:icb2},  the   last statement follows.
\end{proof}
Let us state and prove a special instance of Theorem~\ref{th:cbn}. For that we need a result, which can be found in \cite[Lemma 8.27]{Amir}.
\begin{lemma} \label{le:cbqc}
Let $a>0$, $d\geq 0$ and  $N\geq 1$. Then, 
$$
\displaystyle \frac{d+ a\sum_{k=1}^{N} \frac{1}{k}}{\sum_{k=1}^{N} \frac{1}{\sqrt{k+1}}}\leq \frac{4(d+a+a\ln(N))}{\sqrt{N}}.
$$
\end{lemma} 
\begin{remark}
If $(\gamma _{k})_{k\in \mathbb{N}}$ satisfies  (H4), then $(\gamma _{k})_{k\in \mathbb{N}}$ also  satisfies (H3).    The sequence  $(\gamma^k)_{k\in\mathbb{N}}$ with $\gamma_k=1/(k^{1-\theta/2})$ and  $\theta \in (0,1)$ satisfies (H3). Using Lemma~\ref{le:cbqc} we can also prove that  sequence  $(\gamma^k)_{k\in\mathbb{N}}$ with $\gamma_{k}=1/\sqrt{k}$  satisfies (H3). 
\end{remark} 

The proof of the next theorem  follows by combining  inequality \eqref{eq:icb2} of Theorem~\ref{th:cbn} with Lemma~\eqref{le:cbqc}.
\begin{theorem} %\label{th:cbqc}
Let $x^*\in \Omega^*$,  $(x_{k})_{k\in \mathbb{N}}$ be generated by Algorithm~\ref{Alg1s} with $\rho >1/2$ and 
$$
\gamma_{k}=\frac{1}{\sqrt{k}}, \qquad \forall k\in {\mathbb N}.
$$
Then, for each fixed  $N\in {\mathbb N}$, the following inequality  hold: 
\begin{equation*} 
\min \big\{f(x_k)-f^*: ~k=1, \ldots, N\big\} \leq \frac{4}{\Gamma}\frac{\|x_1-x^*\|^2+\beta \rho^{-1}c +\beta \rho^{-1}c\ln(N)}{\sqrt{N}}.
\end{equation*} 
Consequently,  $\lim_{N\to +\infty}\min \big\{f(x_k)-f^*: ~k=1, \ldots, N\big\}=0$.
\end{theorem} 
\begin{theorem} %\label{th:cb}
Assume that  $\Omega^*\neq \varnothing$. Let $(x_{k})_{k\in \mathbb{N}}$ be generated by Algorithm~\ref{Alg1s} with $\rho>1/2$ and $x^*\in \Omega^*$. Then, for each fixed  $N\in {\mathbb N}$, the following inequality  holds:
\begin{equation} \label{eq:icb1}
\min \big\{f(x_k)-f^*: ~k=1, \ldots, N\big\} \leq \frac{1}{\Gamma} \Big(\|x_1-x^*\|^2+\beta \rho^{-1}c\sum_{j=1}^N\gamma_k^2\Big)\frac{1}{N\gamma_{N+1}}. 
\end{equation} 
As a consequence,  if (H4) holds then $\lim_{N\to +\infty} \min \big\{f(x_k)-f^*: ~k=1, \ldots, N\big\}=0$.
\end{theorem} 
\begin{proof}
Let $k\leq N$. Applying the  inequality \eqref{eq:sfiaux} in  Lemma~\ref{le:auxineq}  we obtain that 
$$
\Gamma \sum_{k=1}^N \left(\gamma_{k+1}(f(x_{k})-f^*)\right)\leq \sum_{k=1}^N\left(\|x_{k}-x^*\|^{2} -\|x_{k+1}-x^*\|^{2}\right)+  \frac{1}{\rho}\beta c \sum_{k=1}^N\gamma_k^2.
$$
Thus,   using (H3) we conclude that 
$$
\Gamma \gamma_{N+1} \sum_{k=1}^N(f(x_{k})-f^*)\leq \|x_{1}-x^*\|^{2} +    \frac{1}{\rho}\beta c \sum_{k=1}^N\gamma_k^2.
$$
Since $N \min \big\{f(x_k)-f^*: ~k=1, \ldots, N\big\} \leq  \sum_{k=1}^N(f(x_{k})-f^*) $ the  inequality \eqref{eq:icb1} follows.  The proof of the last statement of the theorem follows by combining   \eqref{eq:icb1}  with (H4).
\end{proof}
We end this section by showing that  $(x_{k})_{k\in \mathbb{N}}$ generated by Algorithm~\ref{Alg1s} converges to a solution of the problem~\eqref{eq:pr_linear} whenever  $\Omega^*\neq \varnothing$.  To this end, we assume that the sequence  $(\gamma _{k})_{k\in \mathbb{N}}$ satisfies the following conditions:
\begin{enumerate}
\item[(H5)] $\sum_{k=1}^{+\infty}  \gamma_{k}^{2}\leq+\infty$;
\item[(H6)] $\sum_{k=1}^{+\infty}  \gamma_{k}=+\infty$.
\end{enumerate}
\begin{remark}
If $(\gamma _{k})_{k\in \mathbb{N}}$ satisfies  (H5) and (H6), then $(\gamma _{k})_{k\in \mathbb{N}}$ also  satisfies (H3). 
\end{remark}
\begin{theorem} %\label{th:fcrss}
 Let $(x_{k})_{k\in \mathbb{N}}$ be generated by Algorithm~\ref{Alg1s} with $\rho>1/2$.  Assume that  (H5) holds.  If $\Omega^*\neq \varnothing$, then $(x_{k})_{k\in \mathbb{N}}$ is bounded. Moreover, if (H6) hold, then $(x_{k})_{k\in \mathbb{N}}$ converges to a solution of problem~\eqref{eq:pr_linear}.\end{theorem}
\begin{proof}
Let $x \in \Omega^*$.  Using \eqref {eq:sfiaux} in Lemma~\ref{le:auxineq}    we obtain after some algebraic manipulations that 
\begin{equation*}  
\|x_{k+1}-x\|^{2} \leq \|x_{k}-x\|^{2} -\Gamma \gamma_{k+1}(f(x_{k})-f^*) + \frac{1}{\rho}\beta c\gamma_k^2,  \quad \forall k \in \mathbb{N}.
\end{equation*}
Thus, considering that   $f(x_{k})-f^*\geq 0$, for all $k \in \mathbb{N}$, it follows from the last inequality that 
\begin{equation*} 
\|x_{k+1}-x\|^{2}\leq \|x_{k}-x\|^{2}+ \frac{1}{\rho}\beta c\gamma_k^2, \quad \forall k \in \mathbb{N}.
\end{equation*}
Hence,  (H5) together with Definition~\ref {def:fejer} implies that the  sequence $(x_{k})_{k\in \mathbb{N}}$ is  quasi-Fej\'er convergent to $\Omega^*$. Since  $\Omega^*\neq \varnothing$, the item $(i)$ of Lemma~\ref{le:fejer}  implies that $(x_{k})_{k\in \mathbb{N}}$  is bounded and the first statement  is proved. To proceed, define a subsequence  $(x_{k_N})_{N\in \mathbb{N}}$ of the sequence  $(x_{k})_{k\in \mathbb{N}}$ such that 
$$
f(x_{k_N})-f^*:=\min \big\{f(x_k)-f^*: ~k=1, \ldots, N\big\}, \qquad N\in {\mathbb N}.
$$
Since $(x_{k})_{k\in \mathbb{N}}$  is bounded, we conclude that $(x_{k_N})_{N\in \mathbb{N}}$ is also bounded. Without loss of generality we can assume that  $(x_{k_N})_{N\in \mathbb{N}}$  converges. Set ${\bar x}=\lim_{N\to \infty}x_{k_N}$. Under the assumptions (H3) and (H6) we  have  from the last part of Theorem~\ref{th:cbn} that $0=\lim_{N\to + \infty}(f(x_{k_N})-f^*)$. Thus,  using that ${\bar x}=\lim_{N\to \infty}x_{k_N}$,  we conclude that $f({\bar x})=f^*$, which implies that  ${\bar x}\in \Omega^*$. Therefore,  due $(x_{k})_{k\in \mathbb{N}}$ be   quasi-Fej\'er convergent to $\Omega^*$, by applying item~$(ii)$ of Lemma~\ref{le:fejer} we obtain the $(x_{k})_{k\in \mathbb{N}}$ converges to ${\bar x}$, which completes the proof.
\end{proof}
%%%%%%%%%%%%%%%%%%%%%%%%%%%%%%%%
\subsection{Convergence analysis for compact constraint set} %\label{eq:csconvsec}
The aim of  this section  is to  analyze  the behavior of   $(x _{k})_{k\in \mathbb{N}}$ under assumptions  (H1), (H2) (H3) and  one  new additional assumption. The new assumption is  as follows:
\begin{enumerate}
\item[(H7)]  The set ${\cal C}$ is compact.
\end{enumerate}
To state the next  theorem let us  introduce  the following  auxiliary  positive constant
\begin{equation*} 
D\geq \max_{x, y\in {\cal C}}\|x-y\|^2, 
\end{equation*} 
and to prove it  we also need an additional result, which can be found in \cite[Lemma 8.27]{Amir}.
\begin{lemma} \label{le:cscbqc}
Let $a>0$, $d\geq 0$ and  $N\geq 2$. Then, 
$$
\displaystyle \frac{d+ a\sum_{k=\lceil N/2 \rceil}^{N}\frac{1}{k}}{\sum_{k=\lceil N/2 \rceil}^{N} \frac{1}{\sqrt{k+1}}}\leq \frac{4(d+a\ln(3))}{\sqrt{N+2}}.
$$
\end{lemma} 
In the next theorem we show that for suitable choice of the sequence   $(\gamma _{k})_{k\in \mathbb{N}}$  the rate of convergence of Algorithm~\ref{Alg1s}  is ${\cal O}(1/\sqrt{k})$.
\begin{theorem} 
 Let $x^*\in \Omega^*$,  $(x_{k})_{k\in \mathbb{N}}$ be generated by Algorithm~\ref{Alg1s} with $\rho >1/2$ and 
\begin{equation} \label{eq:dgkcs}
\gamma_{k}=\frac{1}{\sqrt{k}}, \qquad \forall k\in {\mathbb N}.
\end{equation} 
Then, for each fixed  $N\in {\mathbb N}$, the following inequality  hold: 
\begin{equation*} 
\min \big\{f(x_k)-f^*: ~k=1, \ldots, N\big\} \leq \frac{4\big(D+    \frac{{\beta c}}{\rho} \ln(3)\big)}{\Gamma \sqrt{N+2}}.
\end{equation*} 
Consequently,  $\lim_{N\to +\infty}\min \big\{f(x_k)-f^*: ~k=1, \ldots, N\big\}=0$.
\end{theorem} 
\begin{proof}
It  follows from \eqref{le:auxineq}  in Lemma~\ref{eq:sfiaux} and definition of $\gamma_{k}$ in  \eqref{eq:dgkcs} that 
\begin{equation*}  
\Gamma \frac{1}{\sqrt{k+1}}(f(x_{k})-f^*)\leq \|x_{k}-x^*\|^{2} -\|x_{k+1}-x^*\|^{2}+  \frac{{\beta c}}{\rho}\frac{1}{k},  \quad \forall k \in \mathbb{N}.
\end{equation*}
Thus, summing this  inequality over $k=\lceil N/2 \rceil, \lceil N/2 \rceil +1, \ldots, N$ we conclude that 
\begin{equation*}  
\Gamma \sum_{k=\lceil N/2 \rceil}^{N}\frac{1}{\sqrt{k+1}}(f(x_{k})-f^*)\leq \|x_{\lceil N/2 \rceil}-x^*\|^{2} -\|x_{N+1}-x^*\|^{2}+  \frac{{\beta c}}{\rho}\sum_{k=\lceil N/2 \rceil}^{N}\frac{1}{k}.
\end{equation*}
Since  $\min \big\{f(x_k)-f^*: ~k=1, \ldots, N\big\}  \sum_{k=\lceil N/2 \rceil}^{N} \frac{1}{\sqrt{k+1}}\leq   \sum_{k=\lceil N/2 \rceil}^{N}\frac{1}{\sqrt{k+1}}(f(x_{k})-f^*)$ and considering that $D\geq \max_{x, y\in {\cal C}}\|x-y\|^2$, we obtain
\begin{equation*}  
\Gamma \min \big\{f(x_k)-f^*: ~k=1, \ldots, N\big\} \sum_{k=\lceil N/2 \rceil}^{N}\frac{1}{\sqrt{k+1}}\leq D+  \frac{{\beta c}}{\rho}\sum_{k=\lceil N/2 \rceil}^{N}\frac{1}{k}.
\end{equation*}
The last inequality implies that 
\begin{equation*}  
 \min \big\{f(x_k)-f^*: ~k=1, \ldots, N\big\}  \leq \frac{D+   \frac{{\beta c}}{\rho}\sum_{k=\lceil N/2 \rceil}^{N}\frac{1}{k}}{\Gamma \sum_{k=\lceil N/2 \rceil}^{N}\frac{1}{\sqrt{k+1}}}, 
\end{equation*}
which combined with Lemma~\ref{le:cscbqc} yields the desired inequality. The second statement of theorem  is an immediate consequence of the first one. 
\end{proof}
\begin{theorem} 
Let $f:\mathbb{R}^{n}\to \mathbb{R} $ be a $\sigma$-strongly convex function and $\sigma >0$.  Let $(x_{k})_{k\in \mathbb{N}}$ be generated by Algorithm~\ref{Alg1s} with $\rho >1/2$, 
$$
\gamma_{k}=\frac{2}{\sigma \beta{\Theta}  k}, \qquad \forall k\in {\mathbb N}.
$$
 and $x^*\in \Omega^*$. Then, for each fixed  $N\in {\mathbb N}$, the following inequality  holds:
\begin{equation*} 
\min \big\{f(x_k)-f(x^*): ~k=1, \ldots, N\big\}  \leq    \frac{8\beta c}{\rho \sigma \beta {\Theta} {\Gamma}}\frac{1}{(N+1)}.
\end{equation*} 
As a consequence,   $\lim_{N\to +\infty}\min \big\{f(x_k)-f(x^*): ~k=1, \ldots, N\big\}=0$.
\end{theorem} 
\begin{proof}
Since  $\gamma_{k}=2/(\sigma \beta{\Theta}  k)$  satisfies  (H3),  it follows from  \eqref{eq:psfiaux} in   Lemma~\ref{le:auxineq} and $\sigma >0$ that 
\begin{equation*} 
\frac{{\Gamma}}{\sigma \beta {\Theta}}(f(x_{k})-f(x^*))\leq \frac{1-\sigma \beta {\Theta}\gamma _{k+1}}{\sigma \beta {\Theta} \gamma _{k+1}}  \|x_{k}-x^*\|^{2}-\frac{1}{\sigma \beta {\Theta}\gamma _{k+1}} \|x_{k+1}-x^*\|^{2}+  \frac{\beta c}{\rho \sigma \beta {\Theta}}\frac{\gamma_k^2}{\gamma_{k+1}}.
\end{equation*}
Taking into account that $\gamma_{k+1}=2/(\sigma \beta{\Theta } (k+1))$, the last last inequality becomes 
\begin{equation*} 
\frac{{\Gamma}}{\sigma \beta {\Theta }}(f(x_{k})-f(x^*)) \leq  \frac{k-1}{2}  \|x_{k}-x^*\|^{2}-\frac{k+1}{2} \|x_{k+1}-x^*\|^{2}+  \frac{2\beta c}{\rho (\sigma \beta {\Theta })^2}\frac{2}{k}.
\end{equation*}
Hence, multiplying  the last inequality by $2k$ we obtain that 
\begin{equation*} 
\frac{{2\Gamma}}{\sigma \beta {\Theta }} k (f(x_{k})-f(x^*))\leq k(k-1)  \|x_{k}-x^*\|^{2}-k(k+1) \|x_{k+1}-x^*\|^{2}+  \frac{8\beta c}{\rho (\sigma \beta {\Theta })^2}.
\end{equation*}
Thus, due to $\min \big\{f(x_k)-f(x^*): ~k=1, \ldots, N\big\}  \sum_{k=1}^Nk\leq   \sum_{k=1}^N  k(f(x_{k})-f(x^*))$,   we have
$$
\frac{{2\Gamma}}{\sigma \beta {\Theta }} \min \big\{f(x_k)-f(x^*): ~k=1, \ldots, N\big\}  \sum_{k=1}^N k\leq  -N(N+1) \|x_{N+1}-x^*\|^{2}+  \frac{8\beta c}{\rho (\sigma \beta {\Theta })^2}N.
$$
Therefore, due to $ \sum_{k=1}^N k=N(N+1)/2$, we conclude  that 
$$
{{\Gamma}} \min \big\{f(x_k)-f(x^*): ~k=1, \ldots, N\big\}  \leq    \frac{8\beta c}{\rho \sigma \beta {\Theta }}\frac{1}{(N+1)}, 
$$
which is equivalent  the  desired inequality. The second statement of theorem  is an immediate consequence of the first one. 
\end{proof}
%%%%%%%%%%%%%%%%%%%%%%%%%%%%%
\section{Illustrative numerical experiments} \label{Sec:NumExp}
In this section we present some examples to illustrate the efficiency of the proposed method comparing its performance with other subgradient methods using classical step size rules. It is not our intention to compete with these classical methods or other problem-specific algorithms, but rather to show that a general approach using our method performs remarkably well in a variety of settings. To this end, we consider the same set of constants in all methods and instances. More precisely, we perform Algorithm~\ref{Alg1s} (subgradient method with non-monotone line search) with $c=1$, $\beta=0.9$, $\rho=0.8$, $\alpha_1=0.1$ and $\gamma_{k}=\frac{\zeta}{\sqrt{k}}$, for all $k\geq 1$, for some values of $\zeta$. The other four subgradient methods use different step sizes $\alpha_k$ described in Table~\ref{tablestepsize}. All the methods start from the same initial point ``\texttt{zeros(n,1)}" which means the zero vector in $\mathbb{R}^n$

\begin{table}
\begin{footnotesize}
\begin{center}
\begin{tabular}{ccccc} \hline
Abbr. & Subgradient method & Step size\\ \hline \vspace{0.1cm}
 Constant step & Constant step size & $\displaystyle\alpha_k=0.1$ \\ \vspace{0.2cm}
 Fixed length & fixed step length & $\displaystyle\alpha_k=\frac{0.2}{||g^k||}$ \\ \vspace{0.2cm}
 Nonsum & Non-summable diminishing step & $\displaystyle\alpha_k=\frac{0.1}{\sqrt{k}}$ \\ \vspace{0.1cm}
Sqrsum nonsum & Square summable nut not summable step & $\displaystyle\alpha_k=\frac{0.5}{k}$ \\ \hline
\end{tabular}
\caption{\footnotesize Comparison methods: step sizes in the subgradient method.}
\label{tablestepsize}
\end{center}
\end{footnotesize}
\end{table}

In each case, simple modifications could be made to to improve the performance of our method, but these examples serve to illustrate an implementation of the proposed method and highlight several features. All numerical experiments are implemented in MATLAB R2020b and executed on a personal laptop (Intel Core i7, 2.30 GHz, 8 GB of RAM).

\subsection{Maximum of a finite collections of linear functions}
The experiments of this section are  generated by the class of functions which are point wise maximum of a finite collections of linear functions. These functions are defined as follows:
\begin{equation}\label{maxfunction}
f(x)=\max \{ f_j(x)={a_j}^{\top} x + b_j \, : \, j=1,\ldots,m\},
\end{equation}
where $a_j\in \mathbb{R}^n$ and $b_j\in\mathbb{R}$. In this case, $\partial f(x) = \mbox{conv}\, \{\partial f_i(x) \, : \, f_i(x)=f(x)\}$.  In this example, we consider the vectors $a_j=(a_{j,1},\ldots,a_{j,n})\in \mathbb{R}^n$ and $b_j\in\mathbb{R}$ randomly chosen by ``\texttt{randn}", a build-in MATLAB function which returns normally distributed random numbers.

As mentioned before, all the methods start from the same initial point and they stop if the iterate $k=3000$ is attained. We compare the performance of the methods for different dimensions $n=2$, $n=5$, $n=10$, $n=20$, $n=50$ and $n=100$, where in Algorithm~\ref{Alg1s} we consider the values of $\zeta$ in $\{\gamma_k\}$ as $0.01$, $0.5$, $1.0$, $0.95$, $1.5$ and $3.3$, respectively. The comparison of the methods is done in terms of the difference $f_{best} - f_{min}$, where the value $f_{best}$ stands to the best value of $f(x^k)$ attained and $f_{min}$ denotes the solution of the problem computed by CVX, a package for specifying and solving convex programming; see \cite{GrantBoyd2014,GrantBoyd2008}.

\begin{figure}[h!]
\centering
\subfloat[$n=2$ and $m=10$]{\label{fig1:a}\includegraphics[width=0.3\linewidth]{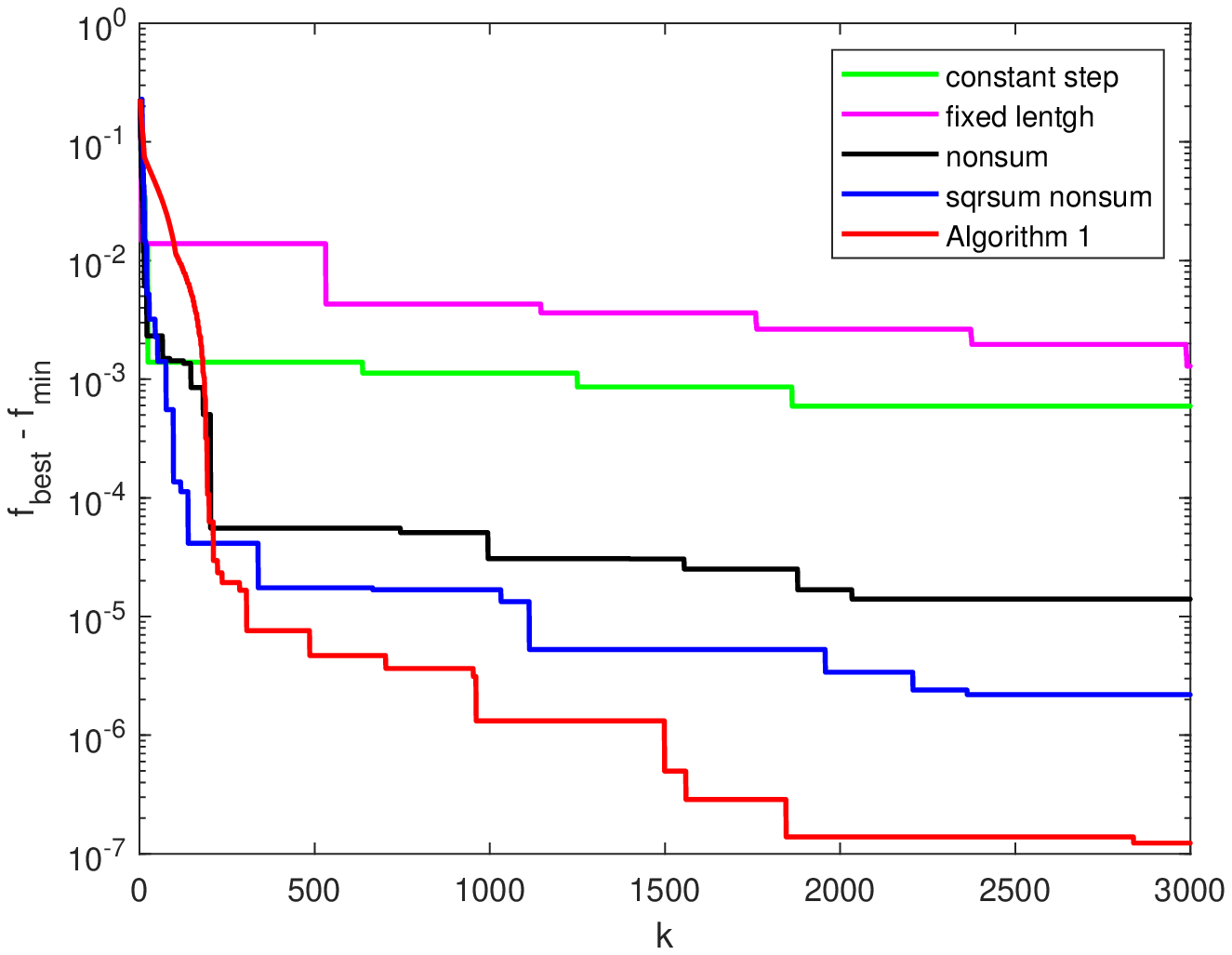}}
\subfloat[$n=5$ and $m=30$]{\label{fig1:b}\includegraphics[width=0.3\linewidth]{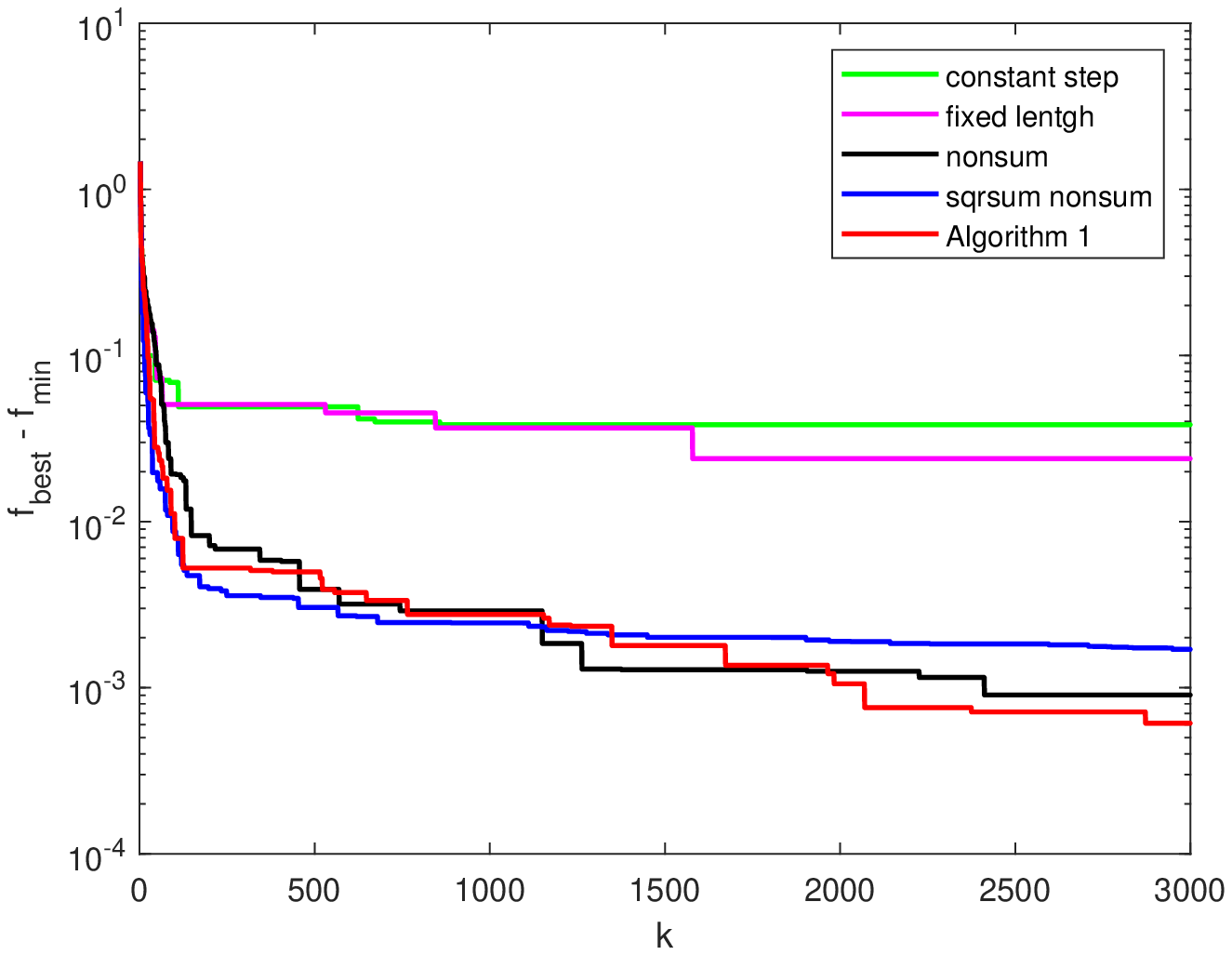}}
\subfloat[$n=10$ and $m=50$]{\label{fig1:c}\includegraphics[width=0.3\linewidth]{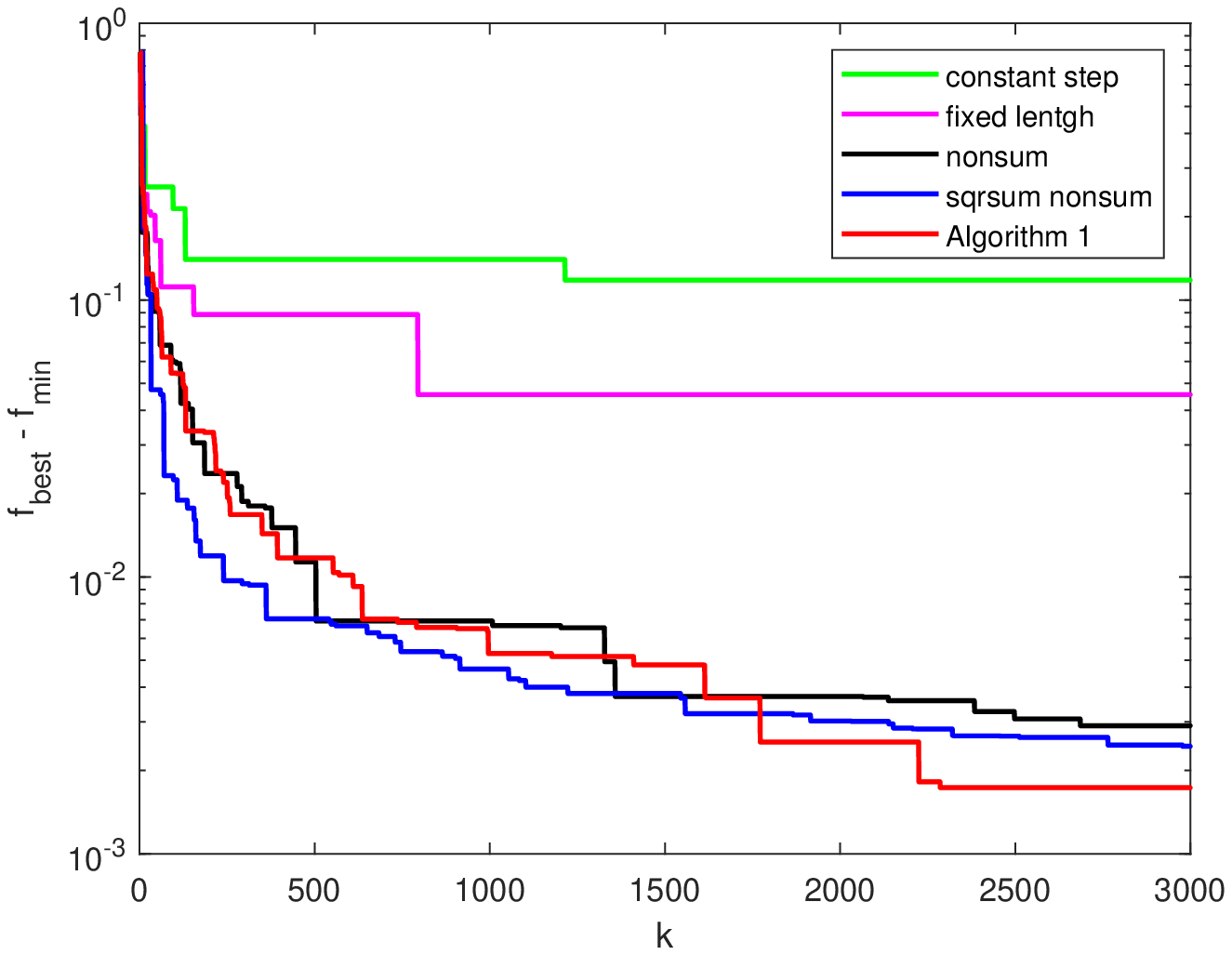}}
\caption{\footnotesize Best value of $f(x^k) - f^*$ (using log. scale) for Algorithm~\ref{Alg1s} and each step size in Table~\ref{tablestepsize}.}
\label{fig1}
\end{figure}

\begin{figure}[h!]
\centering
\subfloat[$n=20$ and $m=100$]{\label{fig2:a}\includegraphics[width=0.3\linewidth]{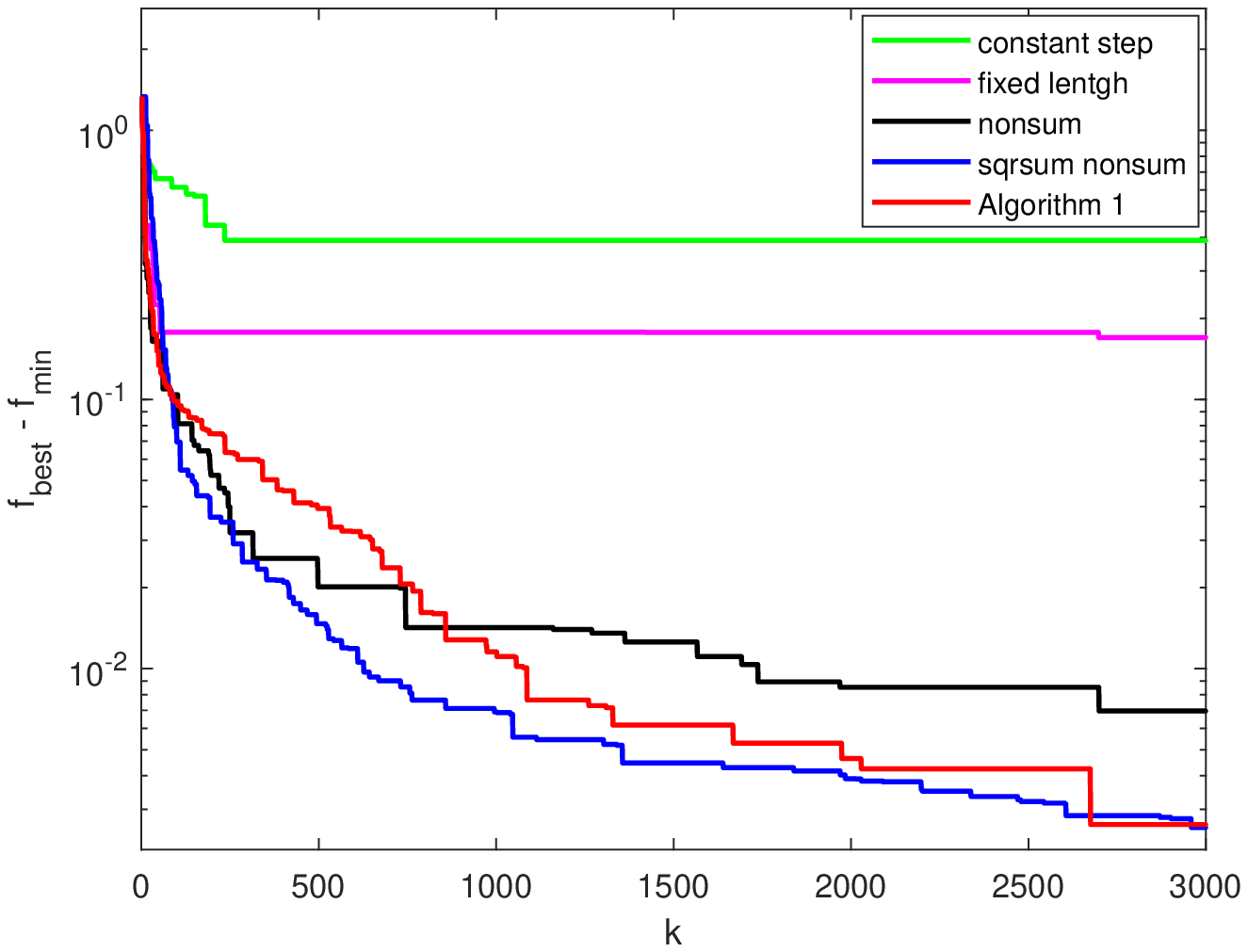}}
\subfloat[$n=50$ and $m=150$]{\label{fig2:b}\includegraphics[width=0.3\linewidth]{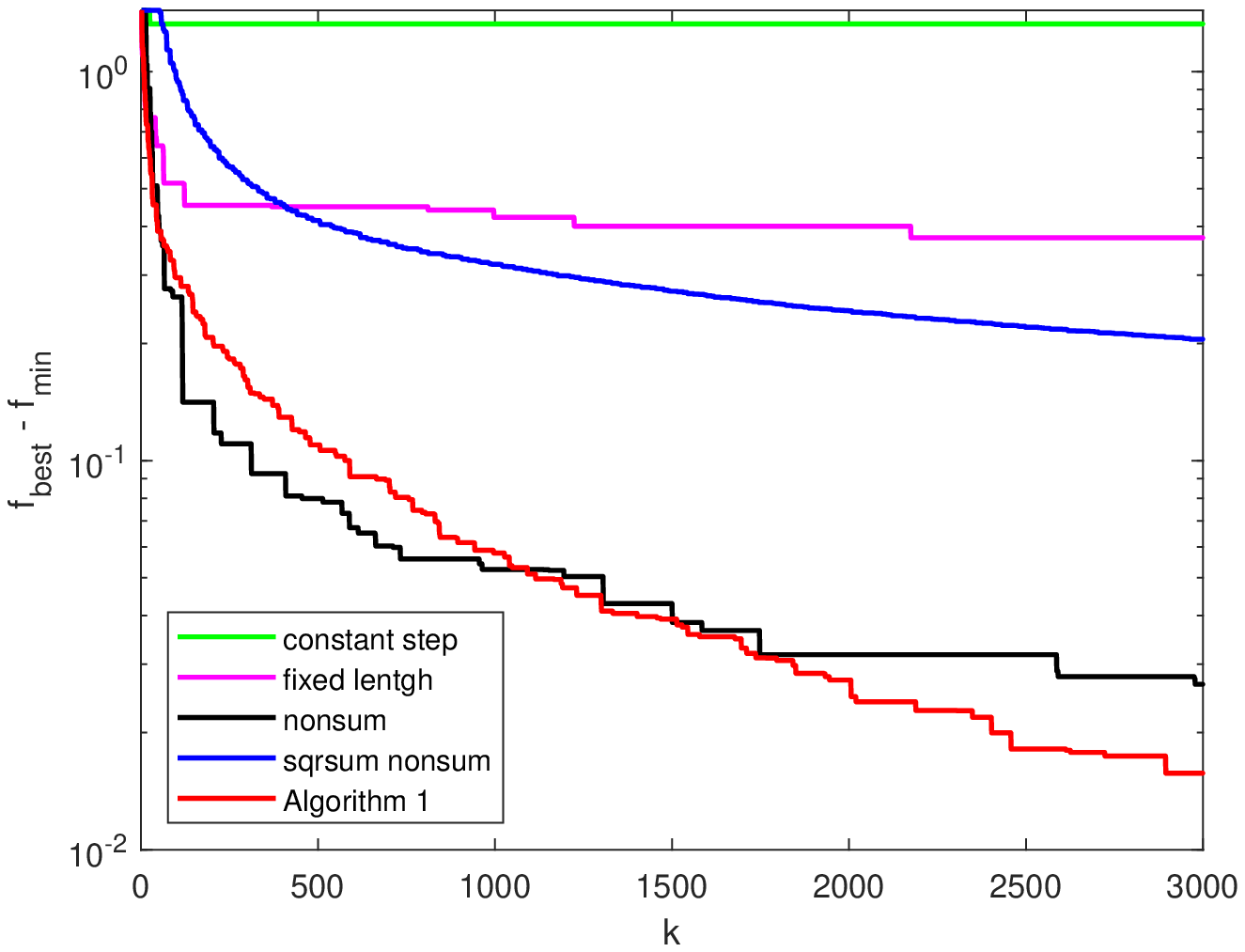}}
\subfloat[$n=100$ and $m=500$]{\label{fig3:c}\includegraphics[width=0.3\linewidth]{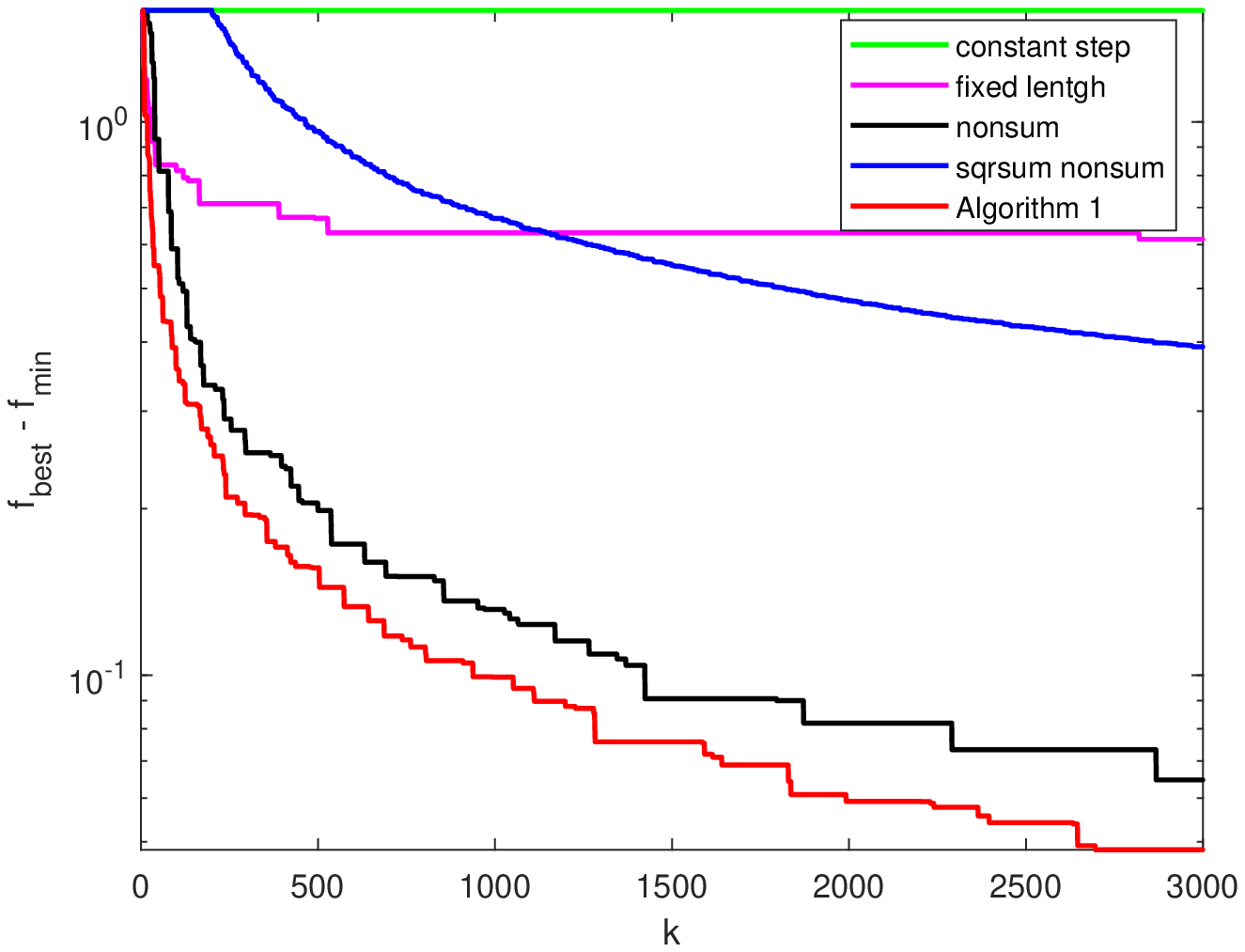}}
\caption{\footnotesize Best value of $f(x^k) - f^*$ (using log. scale) for Algorithm~\ref{Alg1s} and each step size in Table~\ref{tablestepsize}.}
\label{fig2}
\end{figure}

The computation results are displayed in Figures~\ref{fig1}, \ref{fig2} and Tables~\ref{tableresults1},\ref{tableresults2}. In these tables, the first column denotes the dimension $n$ and the number of functions $f_j$, $j=1,\ldots,m$, in \eqref{maxfunction}. The other columns represent, for each method, the best value obtained for $f_{best} - f_{min}$ and the respective iterate $it_{best}$ where it was attained.   As we can see, the results show that Algorithm~\ref{Alg1s} outperforms the other methods providing a better solution or a similar solution in less iterates in all the test problems. In some instances, the subgradient method with the step sizes ``constant step", ``fixed length" and ``sqrsum nonsum" fail to find an acceptable solution in the sense that these methods stop to decrease the objective function in few iterates. In this sense, Algorithm~\ref{Alg1s} and the subgradient method with the step size ``nonsum" have a better performance than the previous ones.

\begin{figure}[h!]
\centering
\subfloat[$n=2$ and $m=10$]{\label{fig3:a}\includegraphics[width=0.3\linewidth]{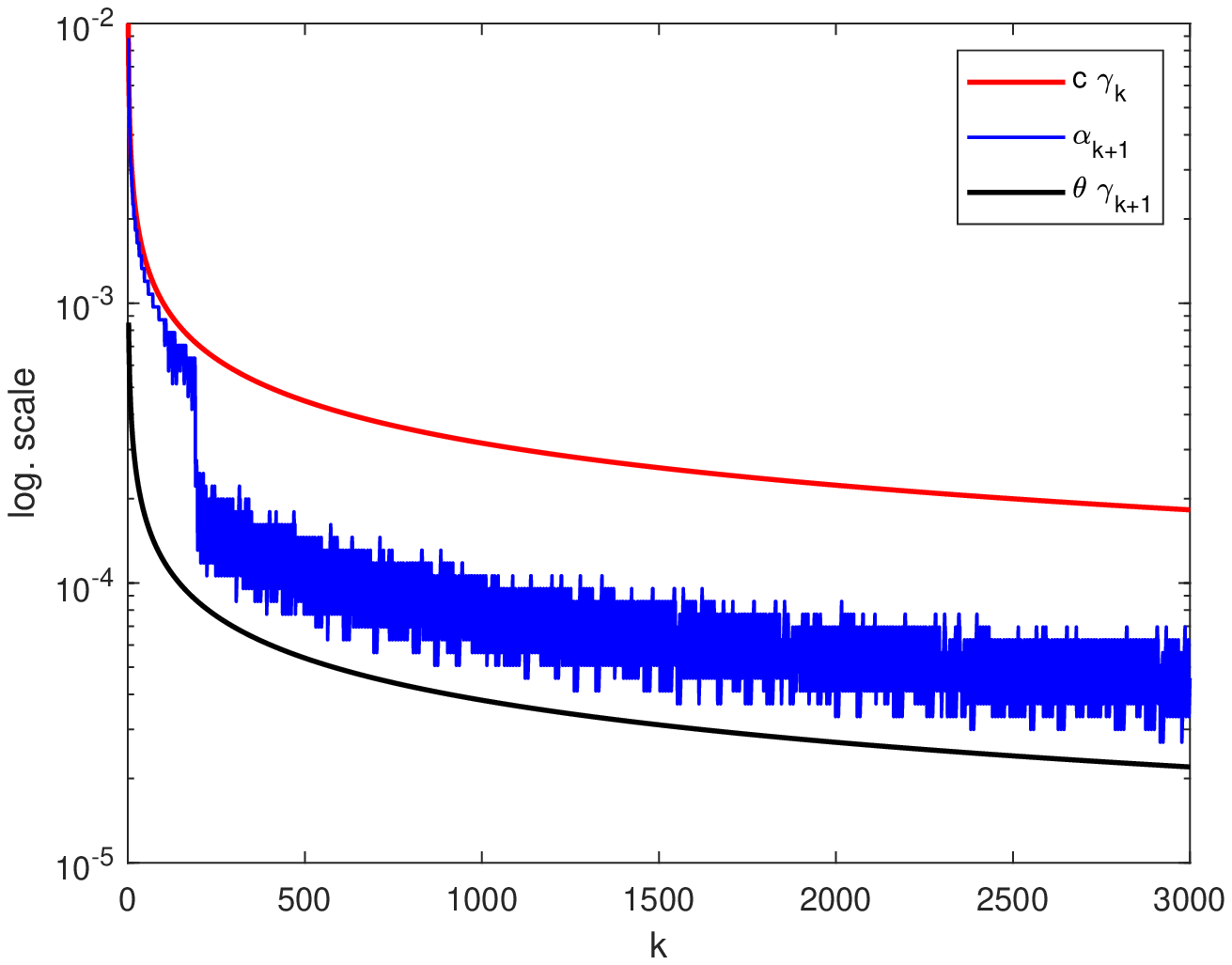}}
\subfloat[$n=5$ and $m=30$]{\label{fig3:b}\includegraphics[width=0.3\linewidth]{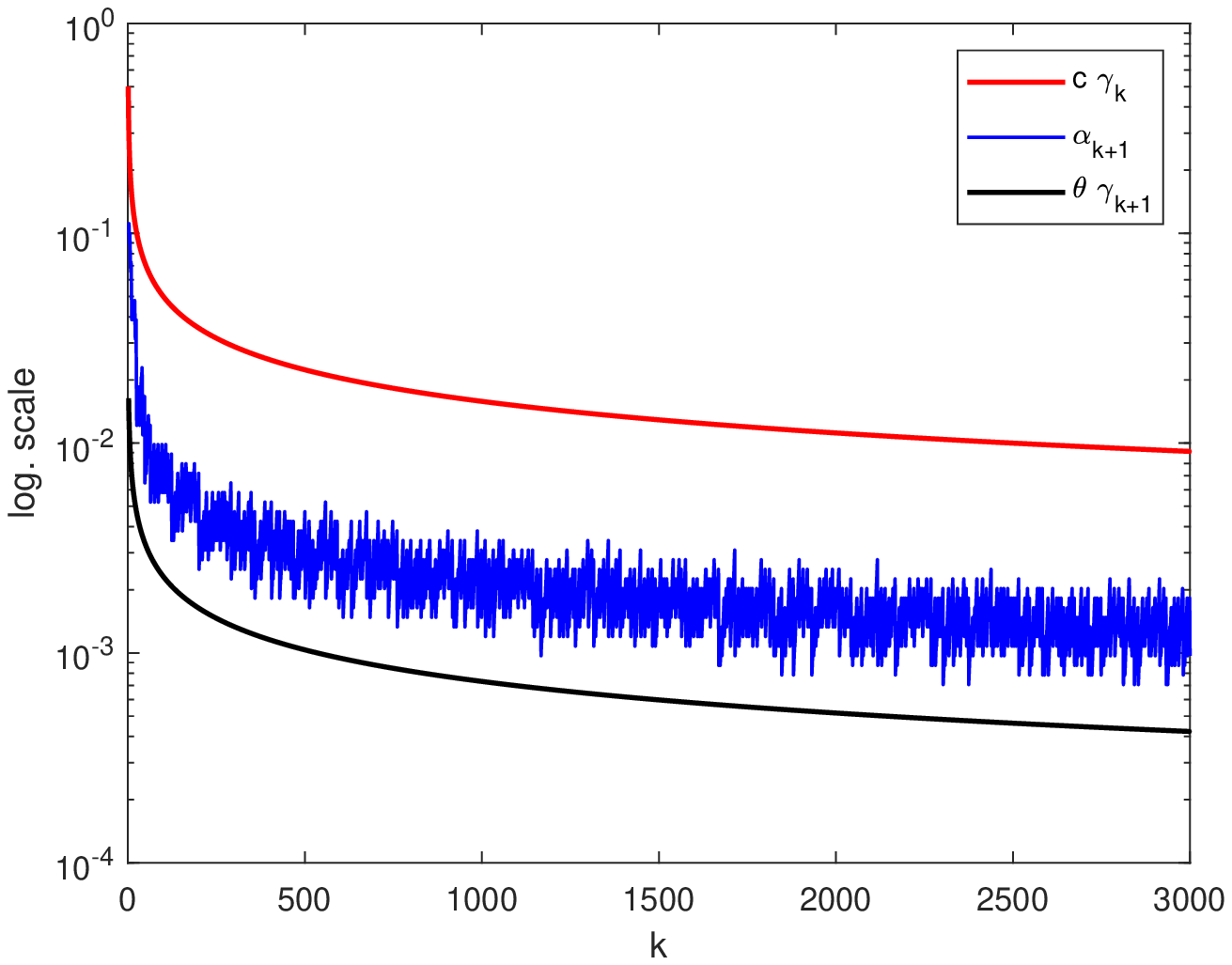}}
\subfloat[$n=10$ and $m=50$]{\label{fig3:c}\includegraphics[width=0.3\linewidth]{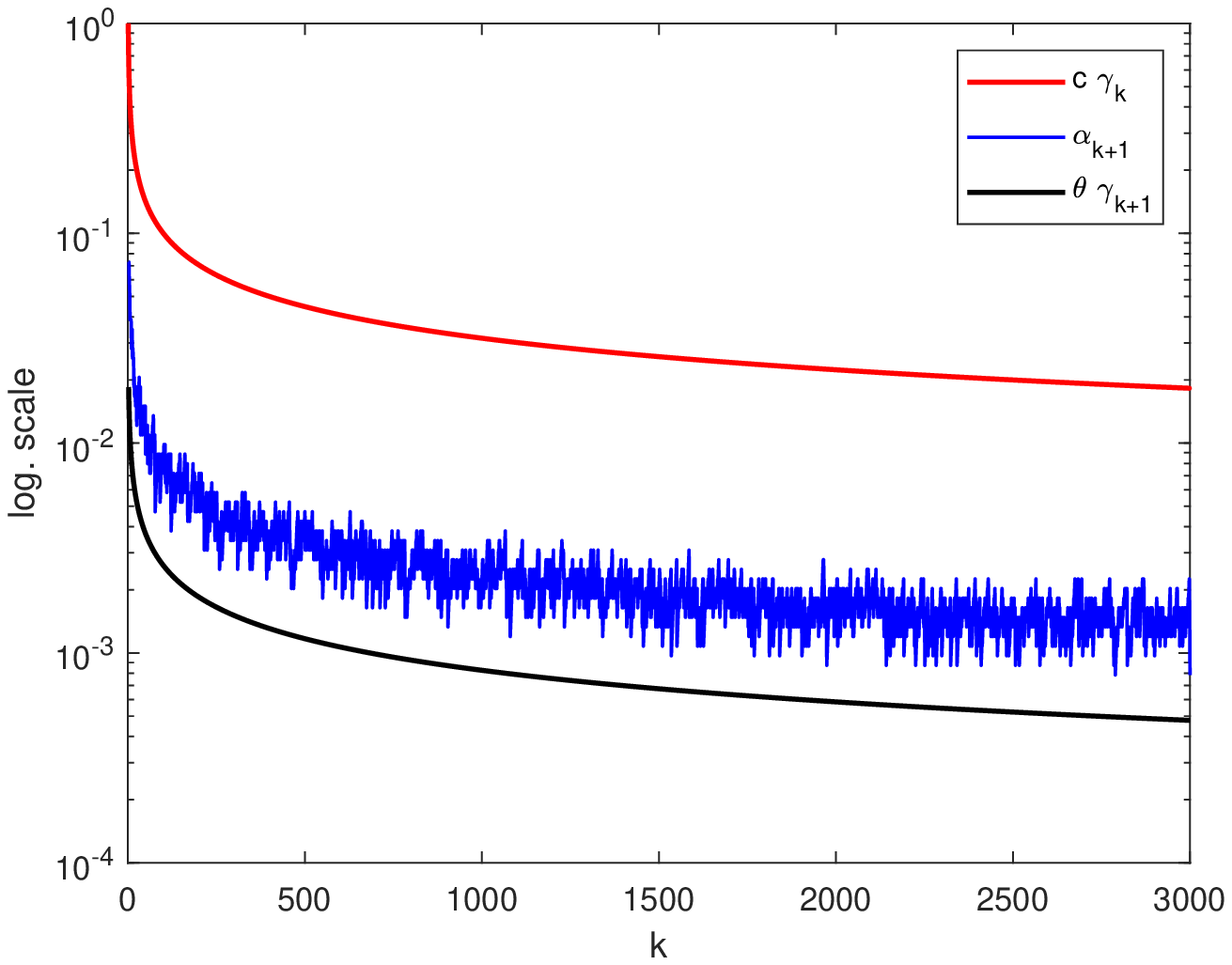}}
\caption{\footnotesize Behavior of the sequences $\{\alpha_k\}$ and $\{\gamma_k\}$ (using log. scale) for Algorithm~\ref{Alg1s}.}
\label{fig3}
\end{figure}

\begin{figure}[h!]
\centering
\subfloat[$n=20$ and $m=100$]{\label{fig4:a}\includegraphics[width=0.3\linewidth]{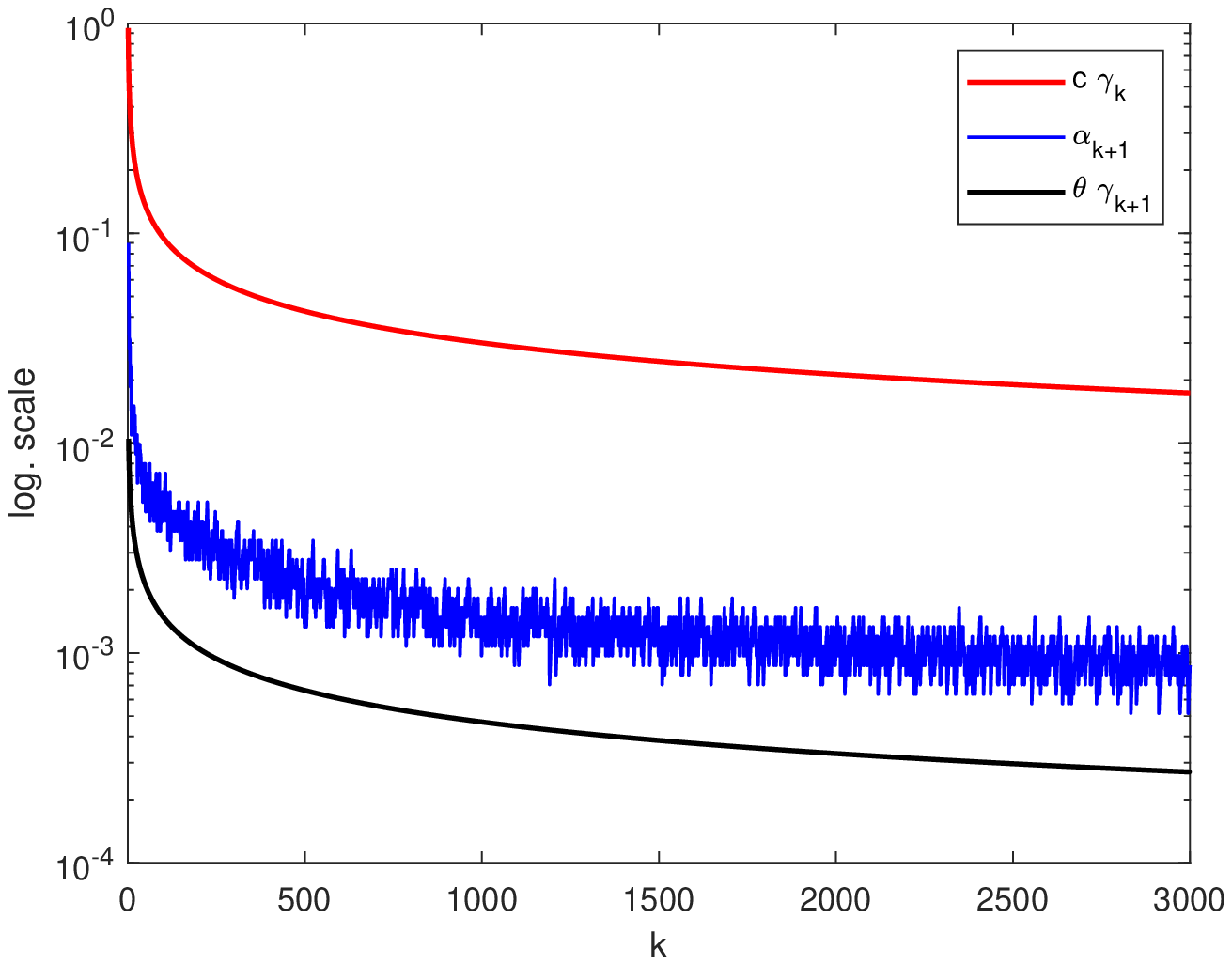}}
\subfloat[$n=50$ and $m=150$]{\label{fig4:b}\includegraphics[width=0.3\linewidth]{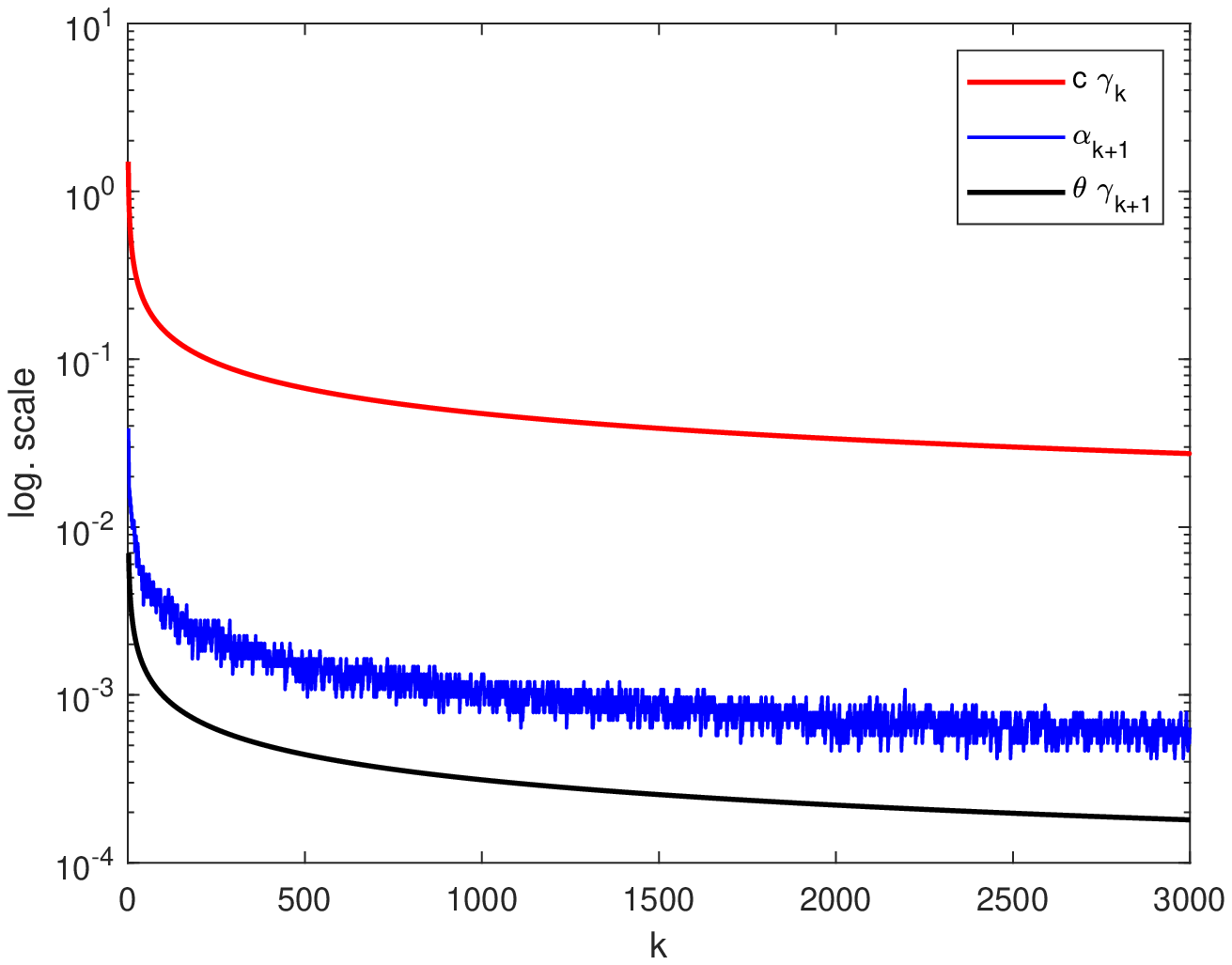}}
\subfloat[$n=100$ and $m=500$]{\label{fig4:c}\includegraphics[width=0.3\linewidth]{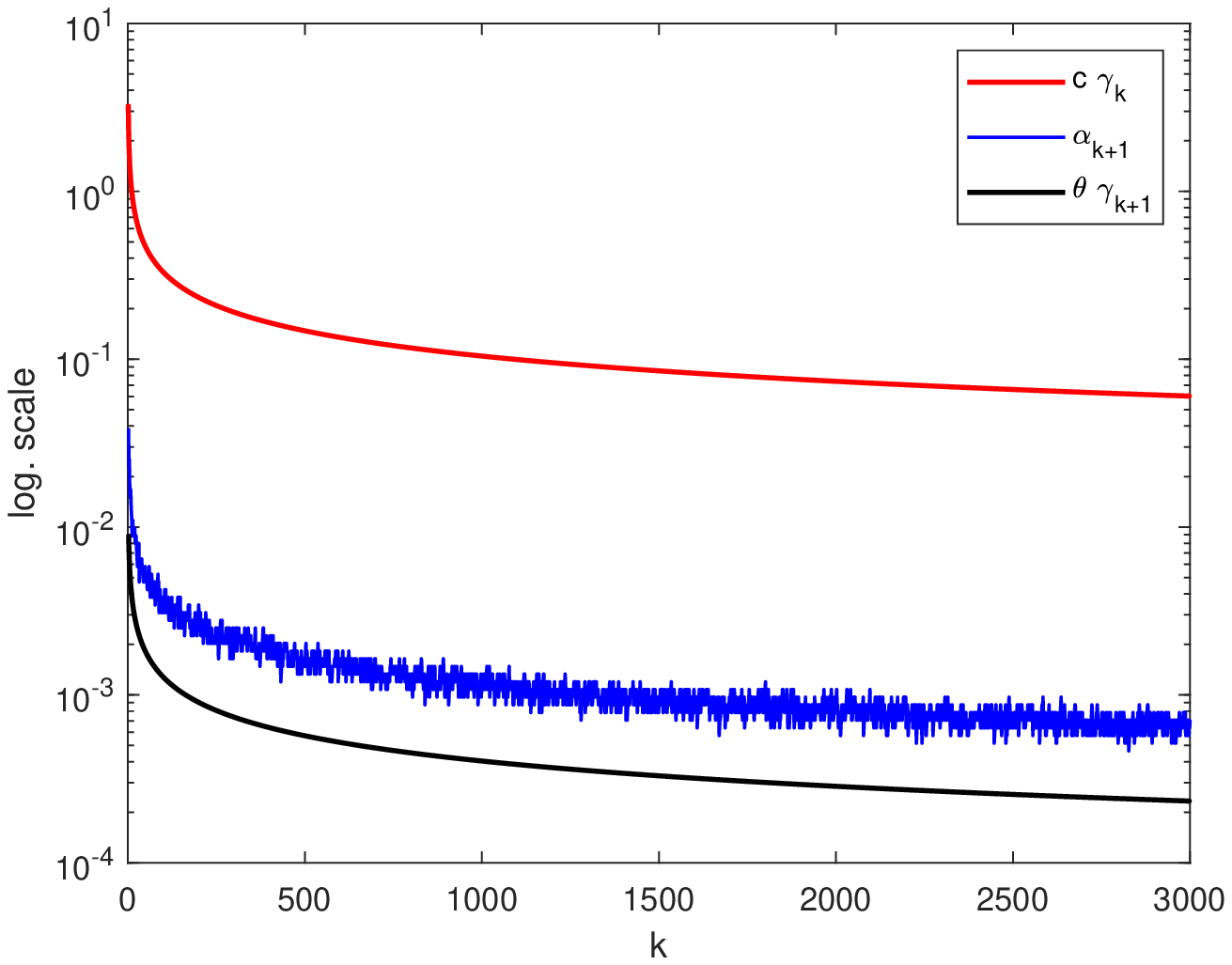}}
\caption{\footnotesize Behavior of the sequences $\{\alpha_k\}$ and $\{\gamma_k\}$ (using log. scale) for Algorithm~\ref{Alg1s}.}
\label{fig4}
\end{figure}

We also investigate the behavior of the sequences $\{\alpha_k\}$ and $\{\gamma_k\}$ in terms of the inequality \ref{eq:bfcnm}, i.e.,
$${\Theta}\gamma_{k+1}\leq \alpha_{k+1}\leq  c\gamma_{k}, \qquad k\in {\mathbb N},$$
where ${\Theta}=\min \left\{ 1,\frac{1}{(1+\rho)L_{f}^{2}} \right\}$. In this example, we consider the Lipschitz constant $L_f$ of the function $f$ in \eqref{maxfunction} as $L_f=\max \{||a_j|| \, : \, j=1,\ldots,m\}$. The results are reported in Figures~\ref{fig3} and \ref{fig4} illustrating the theoretical result stated in Remark~\ref{remarkbfcnm}.

\begin{table}[h!]
\begin{footnotesize}
\begin{center}
\begin{tabular}{l|cc|cc|cc|}
\cline{2-7}
  & \multicolumn{2}{c|}{Algorithm~\ref{Alg1s}} & \multicolumn{2}{c|}{constant step} & \multicolumn{2}{c|}{fixed length}                \\ \cline{2-7} 
 & \multicolumn{1}{c|}{$f_{best} - f_{min}$} & $it_{best}$ & \multicolumn{1}{c|}{$f_{best} - f_{min}$} & $it_{best}$ & \multicolumn{1}{c|}{$f_{best} - f_{min}$} & $it_{best}$  \\ \hline
\multicolumn{1}{|c|}{$n=2$, $m=10$} & \multicolumn{1}{c|}{1.23289e-07}                 &  2838        & \multicolumn{1}{c|}{5.9329e-04}                 &     1863     & \multicolumn{1}{c|}{1.29047e-03}                 &   2990            \\ \hline
\multicolumn{1}{|c|}{$n=5$, $m=30$}          & \multicolumn{1}{c|}{6.11231e-04}                 &  2872        & \multicolumn{1}{c|}{0.038322}                 &    858      & \multicolumn{1}{c|}{0.0239323}                 &     1579      \\ \hline
\multicolumn{1}{|c|}{$n=10$, $m=50$}          & \multicolumn{1}{c|}{1.73369e-03}                 &     2286     & \multicolumn{1}{c|}{0.118048}                 &      1215    & \multicolumn{1}{c|}{0.0455682}                 &   795                 \\ \hline
\multicolumn{1}{|c|}{$n=20$, $m=100$}          & \multicolumn{1}{c|}{2.63594e-03}                 &     2675     & \multicolumn{1}{c|}{0.389903}                 &      235    & \multicolumn{1}{c|}{0.169888}                 &   2698                 \\ \hline
\multicolumn{1}{|c|}{$n=50$, $m=150$}          & \multicolumn{1}{c|}{0.0157351}                 &     2895     & \multicolumn{1}{c|}{1.32626}                 &      25    & \multicolumn{1}{c|}{0.37411}                 &   2175                 \\ \hline
\multicolumn{1}{|c|}{$n=100$, $m=500$}          & \multicolumn{1}{c|}{0.0483826}                 &   2696       & \multicolumn{1}{c|}{1.59047}                 &     2     & \multicolumn{1}{c|}{0.613463}                 &        2820            \\ \hline
\end{tabular}
\caption{\footnotesize Iteration where each algorithm attains the best value of $f(x^k) - f^*$.}
\label{tableresults1}
\end{center}
\end{footnotesize}
\end{table}

\begin{table}[h!]
\begin{footnotesize}
\begin{center}
\begin{tabular}{l|cc|cc|}
\cline{2-5}
  &  \multicolumn{2}{c|}{nonsum} & \multicolumn{2}{c|}{sqrsum nonsum}                \\ \cline{2-5} 
 &  \multicolumn{1}{c|}{$f_{best} - f_{min}$} & $it_{best}$ & \multicolumn{1}{c|}{$f_{best} - f_{min}$} & $it_{best}$  \\ \hline
\multicolumn{1}{|c|}{$n=2$, $m=10$}           & \multicolumn{1}{c|}{1.4022e-05}                 &    2034      & \multicolumn{1}{c|}{2.1927e-06}                 &     2363       \\ \hline
\multicolumn{1}{|c|}{$n=5$, $m=30$}           & \multicolumn{1}{c|}{9.0404e-04}                 &    2412      & \multicolumn{1}{c|}{1.7e-03}                 &     2949      \\ \hline
\multicolumn{1}{|c|}{$n=10$, $m=50$}          &  \multicolumn{1}{c|}{2.90217e-03}                 &      2686    & \multicolumn{1}{c|}{2.4701e-03}                 &        2765            \\ \hline
\multicolumn{1}{|c|}{$n=20$, $m=100$}          &  \multicolumn{1}{c|}{6.96234e-03}                 &      2699    & \multicolumn{1}{c|}{2.57007e-03}                 &        2959            \\ \hline
\multicolumn{1}{|c|}{$n=50$, $m=150$}          &  \multicolumn{1}{c|}{0.0266287}                 &      2978    & \multicolumn{1}{c|}{0.205205}                 &        2964            \\ \hline
\multicolumn{1}{|c|}{$n=100$, $m=500$}          &  \multicolumn{1}{c|}{0.0646978}                 &      2868    & \multicolumn{1}{c|}{0.39232}                 &         2972           \\ \hline
\end{tabular}
\caption{\footnotesize Iteration where each algorithm attains the best value of $f(x^k) - f^*$.}
\label{tableresults2}
\end{center}
\end{footnotesize}
\end{table}

%%%%%%%%%%%%%%%%%%%%%%%%%
\subsection{Fermat-Weber location problem}
The experiment of this section  is the well known Fermat-Weber location problem; see for instance Brimberg~\cite{Brimberg}. Let $a_1, \ldots, a_m$ be given points in $\mathbb{R}^n$. The Fermat-Weber location problem is to solve the following minimization problem
\[
\min_{x\in\mathbb{R}^n}f(x)=\sum_{i=1}^{m}w_i || x - a_i||.
\]
In our particular application, we consider the data points $a_i$, for $i=1, \ldots, 27$, given by the coordinate\footnote{The latitude/longitude coordinates of the Brazilian cities can be found, for instance, at {\it ftp://geoftp.ibge.gov.br/Organizacao/Localidades}.} of the cities which are capital of all 26 states of Brazil and Bras\'ilia (the Federal District, capital of Brazil). We take equally weights for all $a_i$, namely, $w_i=1$, $i=1,\ldots,27$, and consider the integer part of the coordinates converting it from positive to negative to match with the real data. Our goal is to find a point that minimizes the sum of the distances to the given points representing the cities in order to see how distance is such a point from Bras\'ilia (the capital of Brazil).  We denote by 
\begin{equation}\label{fmincvx}
f_{min}=312.9232964118977, \qquad x_{min}=(-45.9630806884547, -12.7465709013343)
\end{equation}
the solution found by the MATLAB package CVX; see \cite{GrantBoyd2014,GrantBoyd2008}. As mentioned in the beginning of this section, we perform Algorithm~\ref{Alg1s} and other four subgradient methods each of them with different step sizes $\alpha_k$ described in Table~\ref{tablestepsize}. All the methods start from the same initial point and they stop at $k=200$ iterates. In Algorithm~\ref{Alg1s}, we take $\zeta=2$ in the definition of $\{\gamma_k\}$.

\begin{figure}[h!]
\centering
\subfloat[Running Algorithm~\ref{Alg1s}.]{\label{fig5:a}\includegraphics[width=0.4\linewidth]{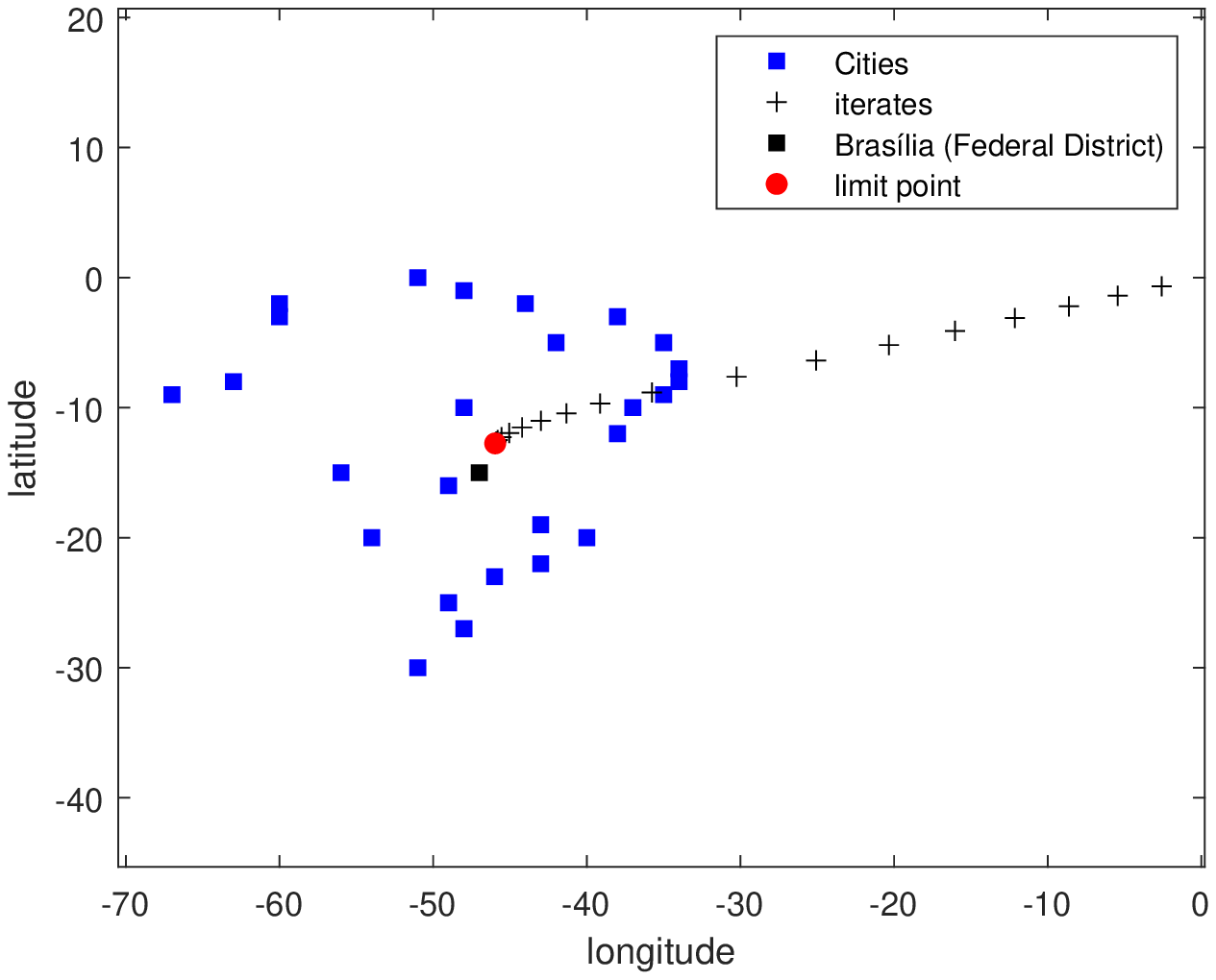}}
\subfloat[Performance of the methods (in log. scale).]{\label{fig5:b}\includegraphics[width=0.4\linewidth]{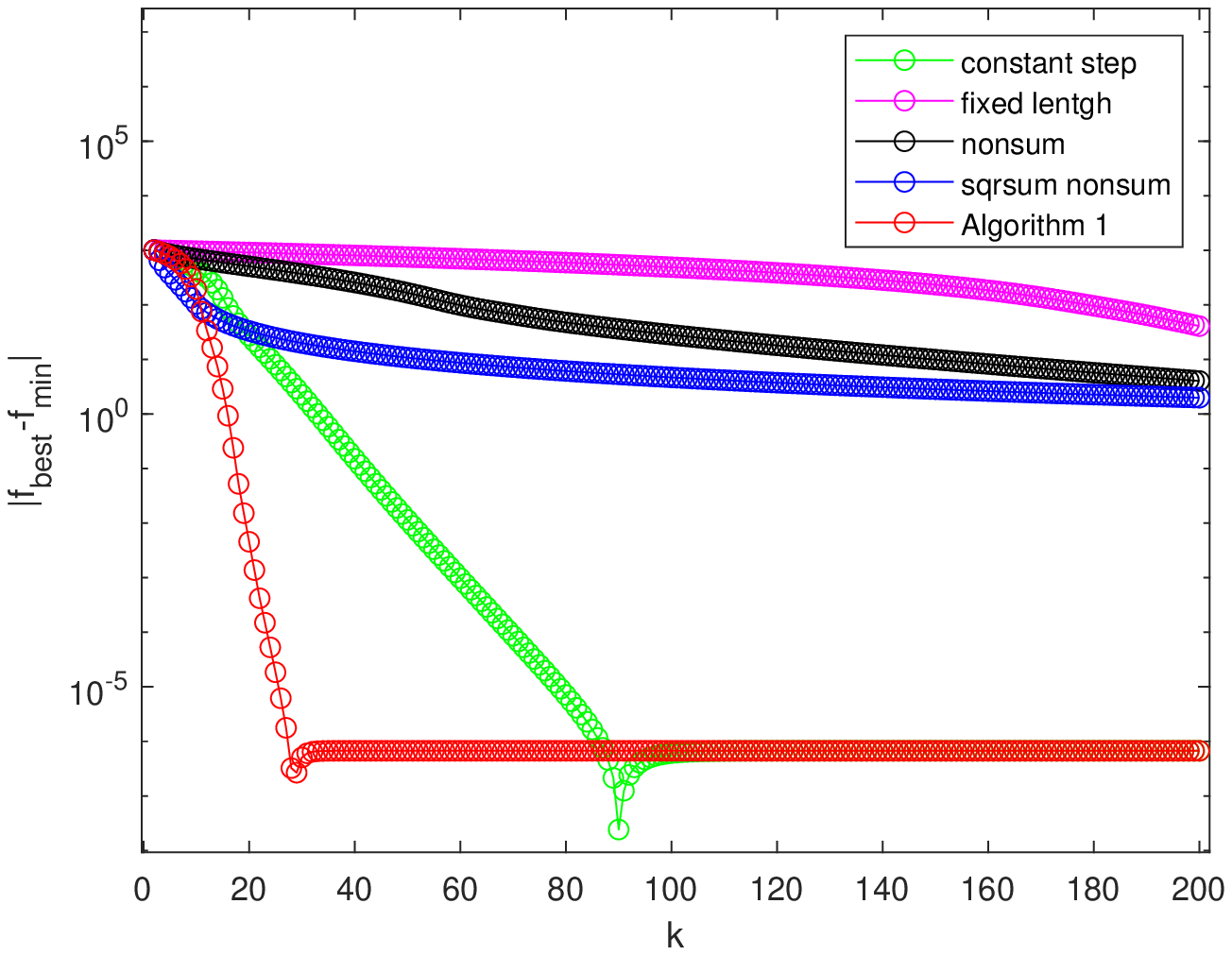}}
\caption{\footnotesize Subgradient methods for solving the Fermat-Weber location problem.}
\label{fig5}
\end{figure}

 In Table~\ref{tableresults3}, we present in the first two columns the solution $x^*=(x_1^*,x_2^*)$ found by each method, in the third column the best value to $|f_{best}-f_{min}|$, where $f_{best}$ stands to the best value of the objective function for each method and $f_{min}$ is given by \eqref{fmincvx}. The last column shows the iterate $it_{best}$ in which the best value $|f_{best}-f_{min}|$ was attained. As we can notice, Algorithm~\ref{Alg1s} and the subgradient method with constant step size found a better solution compared to the solution known ($f_{min}$ and $x_{min}$) than the other methods. However, Algorithm~\ref{Alg1s} found its best value in 29 iterates while the subgradient method with constant step size takes 90 iterates to attain its best value. The performance of each method is presented in Figure~\ref{fig5:b} showing the efficiency of the Algorithm~\ref{Alg1s} for this example. In Figure~\ref{fig5:a}, we present the data of this example as well as the iterates of the Algorithm~\ref{Alg1s} and the solution found by the method.

\begin{table}[h!]
\begin{footnotesize}
\begin{center}
\begin{tabular}{l|cc|cc|}
\cline{2-5}
 % &  \multicolumn{2}{c|}{nonsum} & \multicolumn{2}{c|}{sqrsum nonsum}                \\ \cline{2-5} 
 &  \multicolumn{1}{c|}{$x_1^*$} & $x_2^*$ & \multicolumn{1}{c|}{$|f_{best}-f_{min}|$} & $it_{best}$  \\ \hline
\multicolumn{1}{|c|}{Algorithm~\ref{Alg1s}}           & \multicolumn{1}{c|}{ -45.963064141347097}                 &     -12.746621089909885      & \multicolumn{1}{c|}{2.66879e-07}                 &     29       \\ \hline
\multicolumn{1}{|c|}{constant step}           & \multicolumn{1}{c|}{-45.963064140711523}                 &    -12.746621088320897      & \multicolumn{1}{c|}{2.42824e-08}                 &     90      \\ \hline
\multicolumn{1}{|c|}{fixed length}          &  \multicolumn{1}{c|}{-38.605444422335090}                 &      -9.623064720309808
    & \multicolumn{1}{c|}{40.7379}                 &        200            \\ \hline
\multicolumn{1}{|c|}{Nonsum}          &  \multicolumn{1}{c|}{-43.842367512948982}                 &      -11.429938434104701    & \multicolumn{1}{c|}{ 4.02647}                 &        200            \\ \hline
\multicolumn{1}{|c|}{Sqrsum nonsum}          &  \multicolumn{1}{c|}{-44.521197252917077}                 &       -11.740733447040283    & \multicolumn{1}{c|}{1.9869}                 &            200        \\ \hline
\end{tabular}
\caption{\footnotesize Solution found for the Fermat-Weber location problem.}
\label{tableresults3}
\end{center}
\end{footnotesize}
\end{table}

\section{Conclusions}\label{Sec:Conclusions}
In this paper we have presented a subgradient method with a non-monotone line search for the minimization of convex functions with simple convex constraints. The non-monotone line search allows the method to adaptively select step sizes. As preliminary numerical tests show, this method performs better than the standard subgradient method with prefixed step sizes, which we hope to motivate further research on this subject.

%\section*{Acknowledgments}
%We would like to acknowledge the assistance of volunteers in putting
%together this example manuscript and supplement.

%\bibliographystyle{siamplain}
%\bibliography{bibnmBDCA}

\end{document}